\documentclass[journal]{IEEEtran}

\usepackage{graphicx}
\usepackage{color}
\usepackage{amsfonts}
\usepackage{float}
\usepackage{amsmath}
\usepackage{amsthm}
\usepackage{epstopdf}
\usepackage{changepage}
\usepackage{latexsym}
\usepackage{tikz}
\usepackage{algorithm,algorithmic}
\usepackage{amssymb}
\usepackage{multicol}
\usepackage{cases}
\usepackage{multirow}
\usepackage{tabu}
\usepackage{hyperref}
\usepackage{bm}

\newcommand{\tabincell}[2]{\begin{tabular}{@{}#1@{}}#2\end{tabular}}

\floatstyle{plain}
\newfloat{myalgo}{tbhp}{mya}
{\begin{myalgo}[#1]
\centering
\part{title}\begin{minipage}{#2}
\begin{algorithm}[H]}%
{\end{algorithm}
\end{minipage}
\end{myalgo}}

\begin{document}

\title{Optimal Placement of Dynamic Var Sources by Using Empirical Controllability Covariance}

\author{Junjian~Qi,~\IEEEmembership{Member,~IEEE,}
        Weihong~Huang,~\IEEEmembership{Student Member,~IEEE,}
        Kai~Sun,~\IEEEmembership{Senior Member,~IEEE,}
        and Wei~Kang,~\IEEEmembership{Fellow,~IEEE}
       \thanks{This work was supported in part by U.S. Department of Energy, Office of Electricity Delivery and Energy Reliability under contract DE-AC02-06CH11357, 
       Oak Ridge National Laboratory under the MOVARTI project, NSF CURENT Engineering Research Center (EEC-1041877), and Naval Research Laboratory (Job order RAJH6).
       
    J.~Qi is with the Energy Systems Division, Argonne National Laboratory, Argonne, IL 60439 USA (e-mail: jqi@anl.gov). 
    
    W. Huang and K. Sun are with the Department of Electrical Engineering and Computer Science, University of Tennessee, Knoxville, TN 37996 USA (e-mails: whuang12@utk.edu; kaisun@utk.edu).
    
    W. Kang is with the Department of Applied Mathematics, Naval Postgraduate School, Monterey, CA 93943 USA (e-mail: wkang@nps.edu).

}
}

\maketitle
\markboth{PREPRINT OF DOI: 10.1109/TPWRS.2016.2552481, IEEE TRANSACTIONS ON POWER SYSTEMS.}{stuff}\maketitle

\begin{abstract}
In this paper, the empirical controllability covariance (ECC), which is calculated around the considered operating condition of a power system, 
is applied to quantify the degree of controllability of system voltages under specific dynamic var source locations. 
An optimal dynamic var source placement method addressing fault-induced delayed voltage recovery (FIDVR) issues 
is further formulated as an optimization problem that maximizes the determinant of ECC. 
The optimization problem is effectively solved by the NOMAD solver, which implements the Mesh Adaptive Direct Search algorithm. 
The proposed method is tested on an NPCC 140-bus system and the results show that the proposed method with fault specified ECC can solve the FIDVR issue 
caused by the most severe contingency with fewer dynamic var sources than the Voltage Sensitivity Index (VSI) based method. 
The proposed method with fault unspecified ECC does not depend on the settings of the contingency and can address more FIDVR issues than VSI method 
when placing the same number of SVCs under different fault durations. 
It is also shown that the proposed method can help mitigate voltage collapse. 
\end{abstract}


\section{Introduction}

\IEEEPARstart{T}{here} are increasing interests in optimization of dynamic var sources due to the growing concerns with dynamic voltage security issues, 
especially the fault-induced delayed voltage recovery (FIDVR) issues and even fast voltage collapse at load buses. 
The optimal placement of dynamic var sources such as static var compensators (SVCs) and static synchronous compensator (STATCOM) have been studied in \cite{JGSingh}-\cite{planning}. 

In \cite{JGSingh}, the optimal location of SVCs is determined by a reactive power spot price index under different loading conditions and contingencies. 
 Similarly, the method in \cite{Wibowo} evaluates the annual cost and benefits from FACTS devices installation considering congestion relief and voltage stability. 
In \cite{Minguez}, the placement of SVCs aims at maximizing the loading margin of a transmission network by a multistart Benders decomposition technique. 
In \cite{Farsangi}, the locations of SVCs are selected based on modal analysis and genetic algorithm considering the input signal for supplementary controller of SVCs. 
In [5], the optimal locations and sizing of multi-type FACTS devices are determined by genetic algorithm. 
In \cite{multi}, a multi-contingency constrained reactive power planning method is proposed by decomposing the optimization problem into two phases.
Although these methods can cover different contingencies, they mainly focus on steady-state analysis and cannot capture the full performance of dynamic var sources. 

By contrast, short-term voltage instability is considered in \cite{Vittal}--\cite{planning}. 
In \cite{Vittal}, a trajectory sensitivity index is proposed to identify the location for dynamic var support 
to mitigate short-term voltage instability considering the impact of induction motors for the most severe contingency. 
In \cite{huang} and \cite{huang1}, a voltage sensitivity index is proposed to select the location for reactive power sources to address the same issue in \cite{Vittal}.
In \cite{allocation2}--\cite{Liu}, the linear sensitivity of the performance measure with respect to the reactive power is used to select the candidate location for var sources 
to satisfy the requirements of the voltage stability margin and transient voltage dip under certain contingencies. 
In \cite{allocation1} and \cite{planning}, a sensitivity index is applied to determine the most influential locations of dynamic var support. 
Sensitivity of voltage dip and duration with respect to the addition of var at a specific location is computed along the trajectory of the dynamical system following a disturbance. 
The placements with large overall sensitivity index are chosen as candidate control locations.

Most of the existing methods only consider multiple cases in steady-state analysis or a few contingencies during allocating dynamic var sources in short-term voltage instability issues. Alternatively, the empirical controllability covariance (ECC) \cite{lall}--\cite{hahn}, which has been used in various applications of control system controllability \cite{Krener}, \cite{kang1}, provides a computable tool for empirical analysis of the input-state behavior and make possible the optimal placement of dynamic var sources from the perspective of controllability of nonlinear systems. 

In this paper, we study the placement of dynamic var sources from the perspective of the controllability of nonlinear systems, by applying ECC to quantify 
the degree of controllability of voltage magnitudes by the input from dynamic var sources. 
We formulate the optimal placement of dynamic var sources as an optimization problem that maximizes the determinant of ECC to maximize the controllability of voltages. 

The rest of this paper is organized as follows. Section \ref{s_cont} introduces the fundamentals of controllability and the definition of ECC. 
Section \ref{calECC} presents methods for calculating the ECC used for the dynamic var sources placement. 
Section \ref{allocation} discusses the formulation of optimal placement of dynamic var sources and its implementation. 
Then in Section \ref{vsi method}, we briefly introduce the voltage sensitivity index based placement of dynamic var sources, which 
will be compared with the proposed method in case studies. 
In Section \ref{case}, the proposed method is tested and validated on an NPCC 140-bus system. 
Finally the conclusion is drawn in Section \ref{conclusion}.

\section{Fundamentals of Controllability and Empirical Controllability Covariance} \label{s_cont}

Here, we introduce the fundamentals of controllability and the definition of ECC.

\subsection{Controllability}

For a linear time-invariant system
\begin{subnumcases} {\label{linear}}
\dot{\boldsymbol{x}}=\boldsymbol A\boldsymbol{x} + \boldsymbol{B} \boldsymbol{u} \\
\boldsymbol{y}=\boldsymbol{C}\boldsymbol{x} + \boldsymbol{D}\boldsymbol{u}
\end{subnumcases}
where $\boldsymbol{x} \in \mathbb{R}^n$ is the state vector, $\boldsymbol{u} \in \mathbb{R}^v$ is the input vector, and $\boldsymbol{y}\in \mathbb{R}^p$ is the output vector,
it is controllable if the controllability matrix
\begin{displaymath}
\left[ \begin{array}{ccccc}
\boldsymbol B & \boldsymbol{AB} & \boldsymbol{A}^2 \boldsymbol B & \cdots & \boldsymbol{A}^{n-1} \boldsymbol B
\end{array} \right]
\end{displaymath}
or the controllability gramian \cite{kailath}
\begin{equation}
\boldsymbol{W}_{\textrm{linear}}=\int_0^\infty e^{\boldsymbol{A}t}\boldsymbol{B} \boldsymbol{B}^\top e^{\boldsymbol{A}^\top t}dt
\end{equation}
has full rank.

For a nonlinear system
\begin{subnumcases} {\label{n1}}
\dot{\boldsymbol{x}}=\boldsymbol{f}(\boldsymbol{x},\boldsymbol{u}) \\
\boldsymbol{y}=\boldsymbol{h}(\boldsymbol{x},\boldsymbol{u})
\end{subnumcases}
where $\boldsymbol{x} \in \mathbb{R}^n$ is the state vector, $\boldsymbol{u} \in \mathbb{R}^v$ is the input vector, and $\boldsymbol{y}\in \mathbb{R}^p$ is the output vector,
it is locally controllable at $\boldsymbol{x}_0$ if the nonlinear controllability matrix obtained by
using Lie derivative has full rank at $\boldsymbol{x}=\boldsymbol{x}_0$ \cite{diop}, \cite{diop1}.

The rank test method is easy and straightforward for linear systems. However, for nonlinear systems this can be very complicated even for small systems.
One possibility is to linearize the nonlinear system. But the nonlinear behavior will inevitably be lost.
Alternatively, empirical controllability covariance \cite{lall}, \cite{lall1} provides a computable tool for empirical analysis of the input-state behaviour of a nonlinear system.
It is also proven that the ECC of a stable linear system described by (\ref{linear}) is equal to the usual controllability gramian \cite{lall1}.

\subsection{Empirical Controllability Covariance}  

The following sets are defined for ECC:
\begin{align}
&T^v=\{\boldsymbol{T}_1,\cdots,\boldsymbol{T}_r;\;\;\;\boldsymbol{T}_l \in \mathbb{R}^{v\times v},\;\boldsymbol{T}_l^\top \boldsymbol T_l=\boldsymbol{I}_v,\;l=1,\ldots,r\} \nonumber \\
&M=\{c_1,\cdots,c_s;\;\;\;\;c_m \in \mathbb{R},\;c_m>0,\;m=1,\ldots,s\} \nonumber \\
&E^v=\{\boldsymbol{e}_1,\cdots,\boldsymbol{e}_v;\;\;\;\textrm{standard unit vectors in}\;\mathbb{R}^v\} \nonumber
\end{align}
where $r$ is the number of matrices for excitation directions, $s$ is the number of different excitation sizes for each direction,
$v$ is the number of inputs to the system, and $\boldsymbol{I}_v$ is the identity matrix with dimension $v$.

For a nonlinear system (\ref{n1}), ECC can be defined as \cite{hahn1}
\begin{equation} \label{gd}
\boldsymbol{W}=\sum_{i=1}^{v}\sum_{l=1}^r\sum_{m=1}^{s}\frac{1}{rsc_m^2}\int_0^\infty \boldsymbol{\Psi}^{ilm}(t) dt
\end{equation}
where $\boldsymbol{\Psi}^{ilm}(t)\in \mathbb{R}^{n\times n}$ is given by $\boldsymbol{\Psi}^{ilm}(t)=(\boldsymbol{x}^{ilm}(t)-\boldsymbol{x}_0^{ilm})(\boldsymbol{x}^{ilm}(t)-\boldsymbol{x}_0^{ilm})^\top$, 
$\boldsymbol{x}^{ilm}(t)$ is the state of the nonlinear system corresponding to the input $\boldsymbol{u}(t)=c_m \boldsymbol{T}_l \boldsymbol{e}_i \mathsf{v}(t)+\boldsymbol{u}_0(0)$, and $\mathsf{v}(t)$ is the shape of the input.

For practical implementation, the discrete form of the ECC can be defined as \cite{hahn1}
\begin{equation} \label{gd1}
\boldsymbol{W}=\sum\limits_{i=1}^v\sum\limits_{l=1}^r\sum\limits_{m=1}^s \frac{1}{rsc_m^2}\sum\limits_{k=0}^K \boldsymbol{\Psi}_k^{ilm}\Delta t_k
\end{equation}
where $\boldsymbol{\Psi}_k^{ilm}\in \mathbb{R}^{n\times n}$ is given by $\boldsymbol{\Psi}_k^{ilm}=(\boldsymbol{x}_k^{ilm}-\boldsymbol{x}_0^{ilm})(\boldsymbol{x}_k^{ilm}-\boldsymbol{x}_0^{ilm})^\top$, 
$\boldsymbol{x}_k^{ilm}$ is the state of the nonlinear system at time step $k$ corresponding to the input $\boldsymbol{u}_k=c_m \boldsymbol{T}_l \boldsymbol{e}_i \mathsf{v}_k + \boldsymbol{u}_0(0)$, 
$K$ is the number of points chosen for the approximation of the integral in (\ref{gd}), and $\Delta t_k$ is the time interval between two points.

For multiple inputs $\boldsymbol{u}\in \mathbb{R}^v$, we have
\begin{equation}
\boldsymbol{W}=\sum\limits_{c=1}^v \boldsymbol{W}^c
\end{equation}
where $\boldsymbol{W}^c$ is the empirical covariance for input $c$. 
The ECC for a system with $v$ inputs is the summation of the empirical covariances
computed for each of the $v$ inputs \cite{singh1}.

\section{Calculating Empirical Controllability Covariance for Dynamic Var Sources Placement} \label{calECC}

Here, we discuss how to calculate the ECC used for dynamic var sources placement.  
Because the voltages are of interest for FIDVR issues, we consider the state as the voltage magnitude of all buses in the system, 
which can be obtained by dynamic simulations, such as by using Siemens PTI PSS/E \cite{loadmodel}. 

As mentioned in Section \ref{s_cont}, the ECC for a system with $v$ inputs
is the summation of the empirical covariances for each of the $v$ inputs.
Then the ECC calculated from placing $v$ dynamic var sources individually
adds to be the covariance calculated from placing the $v$ dynamic var sources simultaneously.

Depending on whether or not a fault and the corresponding outputs of the installed dynamic var sources are specified, 
there are two methods to calculate the ECC, which will be discussed separately as follows.

\begin{enumerate} \renewcommand{\labelitemi}{$\bullet$}
\vspace{0.1cm}
\item \textbf{Case 1--fault specified} \\
A specific fault is applied and each time one dynamic var source is installed at one of the load buses, whose reactive power output is considered as the input. 
This is the case in which dynamic var sources need to be placed to address a given fault, e.g. the most severe fault. 
We do not explicitly use the excitation size $c_m$ but install dynamic var sources with different capacities to apply different sizes of input perturbation.
Each input has a series of step changes and at time step $k$ it is:
\begin{equation}
\mathsf{v}_k=Q_{m,k}\,S(t-t_k^{\textrm{step}})
\end{equation}
where $Q_{m,k}$ is the reactive power output of the dynamic var source with the $m$th capacity at time step $k$ and $t_k^{\textrm{step}}$ 
is the time instant when the $k$th step change occurs; $S(t)$ is the unit step function defined as
\begin{displaymath}
S(t) = \left\{ \begin{array}{ll}
0, & \textrm{$t<0,$}\\
1, & \textrm{$t\ge 0.$}
\end{array} \right.
\end{displaymath}

Correspondingly, $c_m$ in (\ref{gd1}) is substituted by $Q_{m,k}$ instead and the definition of ECC is modified to be:
\begin{equation} \label{gd2}
\boldsymbol{W}=\sum\limits_{i=1}^v\sum\limits_{l=1}^r\sum\limits_{m=1}^s \frac{1}{rsQ_{m,k}^2}\sum\limits_{k=0}^K \boldsymbol{\Psi}_k^{ilm}\Delta t_k
\end{equation}
where $T^v=\{\boldsymbol{I}_{v}\}$. 

\vspace{0.1cm}
\item \textbf{Case 2--fault unspecified} \\
No fault is specified and the dynamic response is obtained by perturbing the reactive power load at one load bus. 
The perturbation could be a step, pulse, or some specifically defined function that can well mimic the response of dynamic var sources. 
This generally applies to the case that dynamic var sources need to be placed to address FIDVR issues following multiple faults. 

The case study later will apply the perturbation as illustrated in Fig. \ref{pulse} to potential SVC locations, 
which mimics the var injection from an SVC under fault conditions. 
In Fig. \ref{pulse}, four actual SVC outputs are obtained by simulating the output of a 200 MVar SVC 
installed at four different locations in the NPCC system following the most severe N-1 contingency, 
which can be approximated by a pulse change of 200 MVar lasting for 1s. We may use different magnitudes for the pulse change to calculate the ECC. 

The perturbation can be written as
\begin{equation} \label{pertb}
\mathsf{v}_k=Q_1\, S(t-t_1)-Q_2\, S(t-t_2)
\end{equation}
where $Q_1$, $Q_2$, $t_1$, and $t_2$ are parameters that can be chosen as different values to reflect the response of dynamic var sources under different scenarios, such as different types of faults at different locations with different time durations.

Since the response of dynamic var sources is simulated by changing the reactive power of load buses, for which increasing or reducing reactive power separately correspond to absorbing or injecting reactive power by dynamic var sources, $T^v$ that are used to defined the empirical covariance in (\ref{gd1}) is chosen as $T^v=\{-\boldsymbol{I}_{v}\}$.

\begin{figure}[!t]
\centering
\includegraphics[width=2.8in]{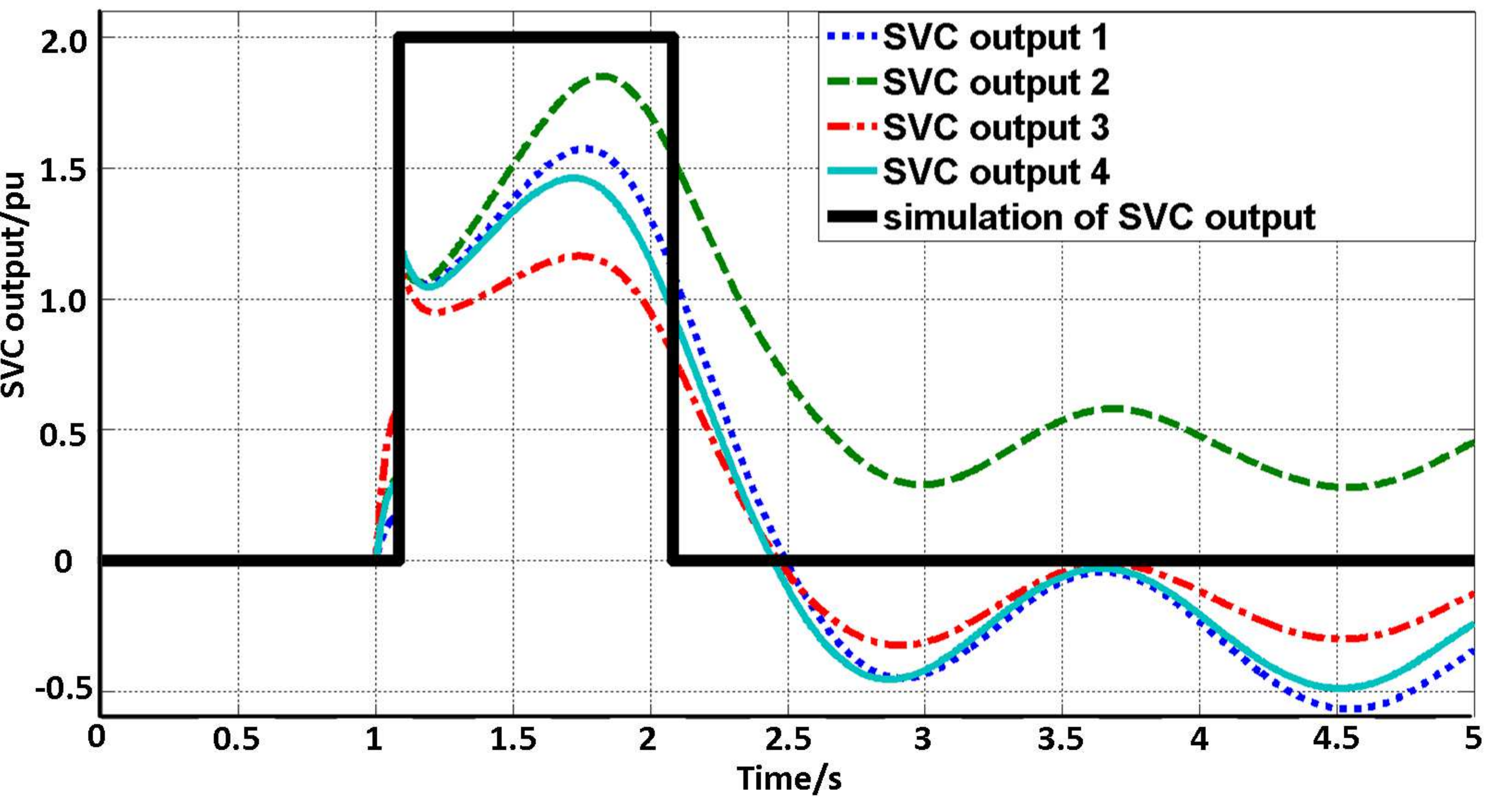}
\caption{Perturbation on reactive power load.}
\label{pulse}
\end{figure}

\end{enumerate}

\section{Dynamic Var Source Placement based on ECC} \label{allocation}

Here, we formulate the optimal dynamic var source placement problem based on ECC and discuss its implementation.

\subsection{Formulation of the Problem}

The degree of controllability can be quantified by a variety of different measures of the ECC, 
such as the smallest eigenvalue \cite{muller}--\cite{qi1}, the trace \cite{singh1}, \cite{qi1}, 
the determinant \cite{muller}, \cite{qi1}--\cite{qi}, or the condition number \cite{Krener}, \cite{qi1}.
Different measures reflect various aspects of controllability: the smallest eigenvalue defines the worst scenario of controllability; 
the trace measures the total gain from input variation to state; the condition number emphasizes the numerical stability; 
and the determinant measures the overall controllability in all directions in noise space.

Although the trace of empirical covariance also tends to measure the overall controllability,
it cannot tell the existence of a zero eigenvalue. Thus an uncontrollable system
may still have a large trace. Compared with the method based on the smallest eigenvalue,
the determinant is a smooth function, which is a desirable property in numerical computations.

Based on these considerations, in this paper we choose the objective as the maximization of the determinant of the ECC under different dynamic var source placements, 
which is the same as that in \cite{qi} for optimal PMU placement for dynamic state estimation.
The optimal placement of dynamic var sources can thus be formulated as
\begin{align} \label{opt}
&\max\limits_{\boldsymbol{z}} \det \, \boldsymbol{W}(\boldsymbol{z}) \nonumber \\
\textrm{s.t.}\; & \;\;\; \sum_{i=1}^L  z_i = v \\
& \;\;\; z_i \in \{0,1\}, \;\; i=1,\ldots,L \nonumber
\end{align}
where $\boldsymbol{z}$ is the vector of binary control variables that determines where to place the dynamic var sources, $\boldsymbol{W}$ is the corresponding ECC, 
$L$ is the number of load buses, and $v$ is the number of dynamic var sources to be placed.

\subsection{Implementation}

Since the ECC for a system with $v$ inputs is the summation of the empirical covariances computed for each of the $v$ inputs individually, as mentioned in Section \ref{s_cont}, 
the ECC calculated from placing $v$ dynamic var sources individually adds to be the identical covariance calculated from placing the $v$ dynamic var sources simultaneously. 
Based on this property, the optimization problem (\ref{opt}) can be written as
\begin{align}\label{opt1}
&\max\limits_{\boldsymbol{z}} \det \, \sum_{i=1}^L z_i \boldsymbol{W}_{i}(\boldsymbol{z}) \nonumber \\
\;\;\textrm{s.t.}\; & \qquad \sum_{i=1}^L  z_i = v  \\
&\qquad z_i \in \{0,1\}, \;\; i=1,\ldots,L \nonumber
\end{align}
where $\boldsymbol{W}_{i}$ is the ECC by only placing one dynamic var source at load bus $i$.

The determinant of a matrix is a high-degree polynomial and its absolute value
can be too small or too huge to be represented as a standard double-precision floating-point number.
By contrast, the logarithm of the determinant is much smaller number and can be simply represented.
Thus we can equivalently rewrite the optimization problem in (\ref{opt1}) as
\begin{align}\label{opt2}
&\min\limits_{\boldsymbol{z}} -\log \det \, \sum_{i=1}^L z_i \boldsymbol{W}_{i}(\boldsymbol{z}) \nonumber \\
\;\;\textrm{s.t.}\; & \qquad \sum_{i=1}^L  z_i = v  \\
&\qquad z_i \in \{0,1\}, \;\; i=1,\ldots,L. \nonumber
\end{align}

To summarize, the optimal placement of dynamic var sources based on the maximization of the determinant of
the ECC can be implemented in the following two steps.

\begin{enumerate} \renewcommand{\labelitemi}{$\bullet$}
\vspace{0.1cm}
\item \textbf{Calculate ECC} \\
The ECC in (\ref{gd}) is calculated by emgr (Empirical Gramian Framework) \cite{emgr} on time interval $[0,t_f]$. We choose $t_f$ as 5 seconds.
Note that we only need to calculate the ECC for placing one dynamic var sources at
one of the load buses and there is no need to compute all the combinations of dynamic var source locations.

\vspace{0.1cm}
\item \textbf{Solve MAX-DET optimization problem} \\
The mixed-integer MAX-DET optimization problem in (\ref{opt2}) is
solved by using the NOMAD solver \cite{nomad},
which is a derivative-free global mixed integer nonlinear programming solver and is
called by the OPTI toolbox \cite{opti}. 
The NOMAD solver implements the Mesh Adaptive Direct Search (MADS) algorithm \cite{mads},
a derivative-free direct search method with a rigorous convergence theory based on the nonsmooth calculus \cite{clarke}, and aims for the best possible solution with a small number of evaluations. Under mild hypotheses, the algorithm globally converges to a point satisfying local optimality conditions based on local properties of the functions defining the problem \cite{nomad}.

NOMAD also includes a Variable Neighborhood Search (VNS) algorithm \cite{audet}, which is based on the VNS metaheuristic \cite{vns}. 
This search strategy perturbs the current iterate and conducts poll-like descents from the perturbed 
point, allowing an escape from local optima on which the algorithm may be trapped \cite{nomad}. 
\end{enumerate}

\section{Voltage Sensitivity Index Based Placement of Dynamic Var Sources} \label{vsi method}

A traditional approach for the placement of dynamic var sources is to evaluate, for a small amount of var injection at a candidate bus, 
how sensitively voltage at one or several target voltage-vulnerable buses will change following the most severe contingency \cite{Vittal}, \cite{huang}. 

The most severe contingency can be identified by simulating a list of credible contingencies, 
which are short-circuit faults near a load bus cleared after a protection reaction time, and 
finding the buses where violations of post-fault voltage criteria happen. 
Some criteria, such as the NERC/WECC reliability standards \cite{standard}, usually limit the post-fault voltage deviations. Examples of such criteria are:

\begin{enumerate} \renewcommand{\labelitemi}{$\bullet$}
\item Post-fault transient voltage dip or overshoot should not exceed 25\% at load buses or 30\% at generator buses, 
and should not exceed 20\% for more than 20 cycles at load buses;
\item Post-transient voltage deviation should not exceed 5\% at any bus.
\end{enumerate}

\begin{figure}[!t]
\centering
\includegraphics[width=2.6in]{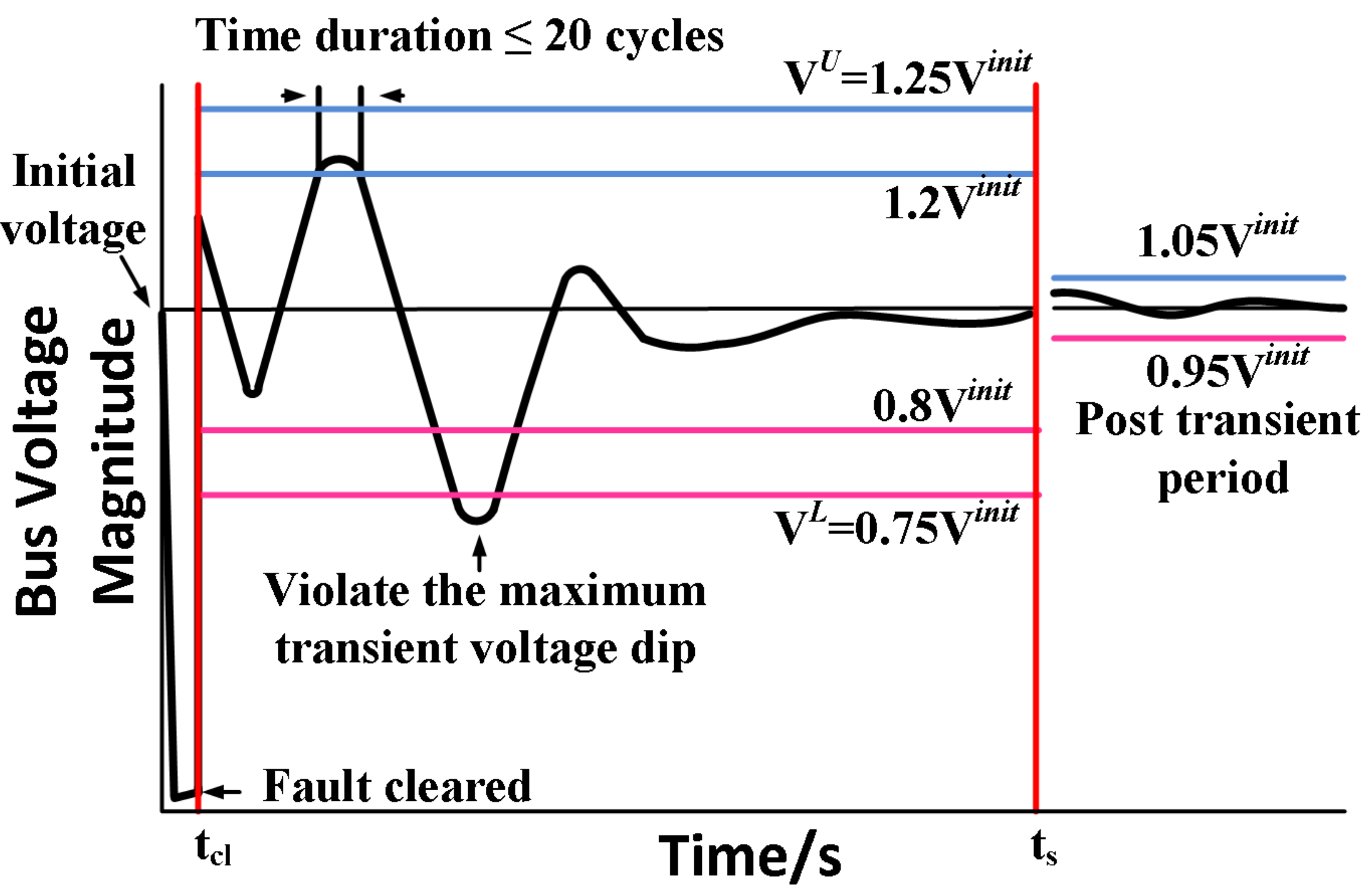}
\caption{Post-fault voltage performance criteria for load bus.}
\label{NERC}
\end{figure}

Fig. \ref{NERC} illustrates the criteria for load buses. Assume that the simulation time step is $\Delta t$, 
the total simulation period is $T \Delta t$. If the $k$th contingency causes violation of the criteria at bus $j$, 
a severity index $\textrm{SI}_{kj}^t$ about the violation at time step $t$ can be evaluated as follows.  

For each bus $j$, the percentage voltage deviation at the $t$th time step is defined as
\begin{equation} \label{Rkj}
R_{kj}^t=\left|\dfrac{V_{kj}^t-V_j^0}{V_j^0}\right| \times 100 \%,\, j=1,\ldots,N,t=1,\ldots,T  
\end{equation}
where $V_{kj}^t$ is the bus $j$ voltage magnitude at time step $t$ and $V_j^0$ is the pre-fault initial voltage magnitude. 
Let $\textrm{SI}_{kj}^t=R_{kj}^t$ if any of the above criteria is violated, or 0, otherwise. 
Then, this average severity index over the simulation period is calculated for each contingency $k$ as
\begin{equation} \label{SIk}
\textrm{SI}_k=\dfrac{1}{T}\sum_{t=1}^T \bigg(\dfrac{1}{N} \sum_{j=1}^{N} \textrm{SI}_{kj}^t \bigg).
\end{equation}

The contingencies with the largest $\textrm{SI}_k$ are considered as the most severe contingencies. 
For an identified severe contingency, a dynamic var source is generally placed 
where var injection causes large sensitivities in terms of post-fault voltage improvements at those buses violating the criteria. 

For instance, a method based on the post-fault Voltage Sensitivity Index (VSI) defined in \cite{huang} may be applied. 
For a severe contingency $k$, a small amount of dynamic var $q_i$ is injected at bus $i$, which can be simulated by adding an SVC or STATCOM of size $q_i$. 
Assume that the var injection changes the voltage at bus $j$ from $V_j^{k,old,t}$ to $V_j^{k,new,t}$ at time step $t$ according to time-domain simulation. 
For each bus $j$ over the whole simulation period, VSI is defined as the maximum voltage recovery sensitivity and can be calculated by
\begin{equation} \label{VSIij}
\textrm{VSI}_{ij}^k=\frac{\max\{V_j^{k,new,t}-V_j^{k,old,t}, t=1,\ldots,T\}}{q_i}.  
\end{equation}
Then, the normalized average voltage sensitivity index for all of the $N$ buses of the system following the injection $q_i$ under contingency $k$ is
\begin{equation} \label{VSIi}
\textrm{VSI}_i^k=\frac{\dfrac{1}{N}\sum\limits_{j=1}^N \textrm{VSI}_{ij}^k}{\max \Big\{\dfrac{1}{N}\sum\limits_{j=1}^N \textrm{VSI}_{lj}^k, l \in \mathcal{C} \Big\}}
\end{equation}
where $\mathcal{C}$ is the set of the candidate buses to install dynamic var sources.

If the top $K$ most severe contingencies are considered, the overall sensitivity index can be calculated as
\begin{equation}
\textrm{VSI}_i=\sum\limits_{k=1}^K \textrm{SI}_k \cdot \textrm{VSI}_i^k.
\end{equation}
If $v$ dynamic var sources are to be installed, using the top $K$ most severe contingency, 
the VSI method calculates $\textrm{VSI}_i$ for each bus and selects the top-$v$ buses with the highest $\textrm{VSI}_i$ to install dynamic var sources.

\section{Case Study} \label{case}

\begin{figure}[!t]
\centering
\includegraphics[width=3.3in]{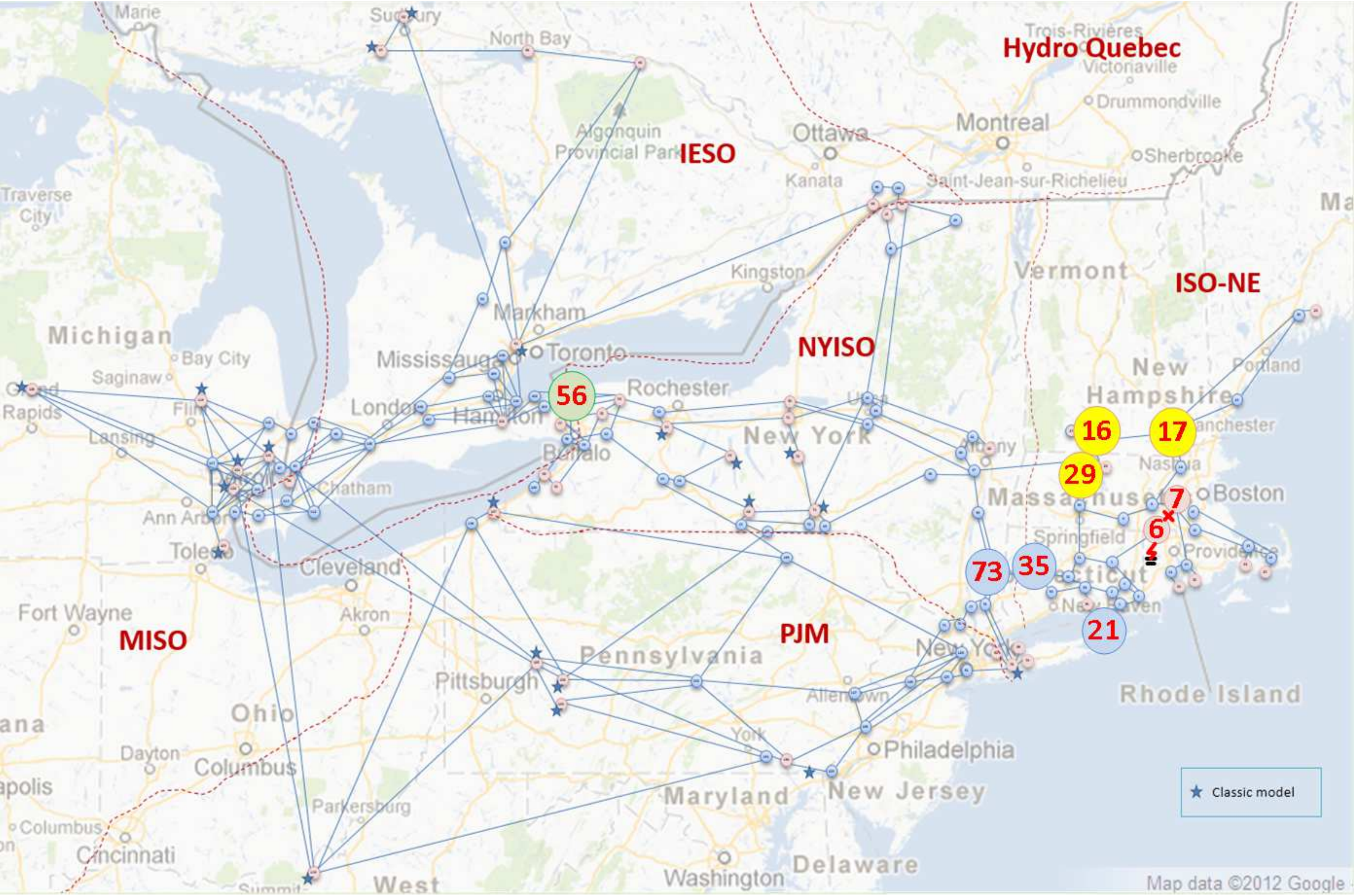}
\caption{NPCC 140-bus system (some important buses are highlighted which will be introduced in the following discussion).}
\label{NPCC1}
\end{figure}

Here, the proposed dynamic var sources placement method is tested on an NPCC 140-bus System \cite{pst}. 
As shown in Fig. \ref{NPCC1}, the NPCC system has 48 generators and 140 buses and represents the northeast region of the Eastern Interconnection system.
For both aforementioned case 1 (fault specified) and case 2 (fault unspecified) in Section \ref{calECC}, the VSI and ECC based methods are compared. 
Note that the SVC placements identified by the VSI-based method can be further improved during the optimization of dynamic var sizes, 
as has been discussed in \cite{allocation1}. Similarly, the SVC locations determined by the proposed ECC-based method can also be improved. 
However, in this paper we do not consider the optimization of dynamic var sizes but only target at finding the optimal placements of SVCs which are assumed to have a predetermined size.

\subsection{Load and SVC Modeling}

We assume that all dynamic var sources are SVCs and all of the SVCs have an identical capacity of 200 Mvar unless otherwise stated. 
In the time-domain simulation, the composite load model CLODBL and SVC model CSVGN5 provided by Siemens PTI PSS/E \cite{loadmodel} 
are applied to simulate FIDVR issues and to provide the solutions by SVCs. 

\begin{figure}[!t]
\centering
\includegraphics[width=3.0in]{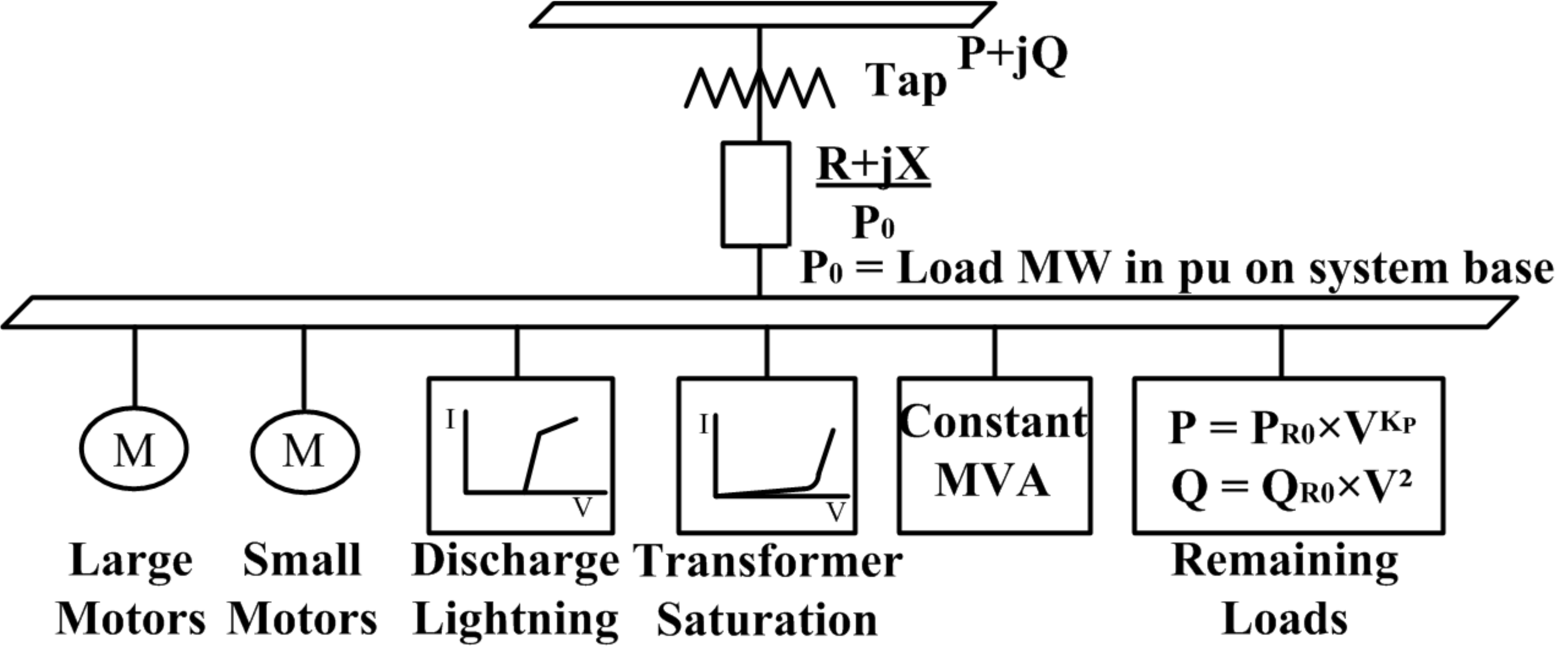}
\caption{Composite load model.}
\label{lm}
\end{figure}

The CLODBL load model shown in Fig. \ref{lm} is used for power system planning and operation studies in PSS/E. 
It consists of induction motors, lighting, and other types of equipment, and is used where it is desirable to represent loads at dynamic level. 
By using this model, users are allowed to specify a minimum amount of data stating the general characteristic of the composite load, 
which can be used internally to establish the relative sizes of motors modeled in dynamic detail. 
The parameters of the composite load are listed in Table \ref{para}.

\begin{table}[H]
\renewcommand{\arraystretch}{1.1}
\caption{Parameters for Composite Load Model}
\label{para}
\centering
\begin{tabular}{cc}
\hline
Parameter & Value \\
\hline
PC\_LM, \% large motor & 0 \\
PC\_SM, \% small motor & 20 \\
PC\_TX, \% transformer exciting current & 1 \\
PC\_DL, \% discharge lighting & 20 \\
PC\_CP, \% constant power & 5 \\
Kp & 1 \\
R, branch resistance in pu & 0 \\
X, branch reactance in pu & 0.1 \\
\hline
\end{tabular}
\end{table}

The CSVGN5 SVC model is a static var system model written for a corresponding model in the WECC stability program in PSS/E. 
It is represented as a generator in power flow simulation. 
In dynamic simulation, CSVGN5 has fast override capability and can be activated when the voltage deviation exceeds a threshold. 
It does not separate the equipment to identify capacitor banks and reactors. 
To maintain the WECC model structure and to include frequency dependence, CSVGN5 assumes that the output is equal to $B_{\textrm{max}}$ and the thyristor-controlled reactor is shut off.
If $B_{\textrm{max}}$ is positive, the capacitor banks are equal to $B_{\textrm{max}}$ times the MVA rating in the power flow from the generator setup. 
If $B_{\textrm{max}}$ is negative, the equipment is assumed to only consist of reactor.  

By using the above models, all N-1 contingencies are simulated. 
For each of them, a three-phase fault is applied to one end of a line and is cleared by opening the line after 5 cycles. 
The total simulation time is 5 seconds and the post transient voltage limits are checked after the first 3 seconds of transient period. Among all simulated contingencies, 
a total of 40 contingencies can cause FIDVR issues.

\subsection{Calculating ECC}

For the fault specified case in Section \ref{calECC}, a specific fault is applied to the system. 
Each time one SVC is installed at one of the candidate load buses. 
Since the typical size of SVCs installed in the transmission level is 200 Mvar, the maximum capacity is selected as 200 Mvar. 
In order to accurately capture the input-state behavior of the system, the size of SVC is selected as a geometric sequence starting from 10 Mvar with a maximum capacity of 
200 Mvar, which are specifically 10 Mvar, 20 Mvar, 40 Mvar, 80 Mvar, 160 Mvar, and 200 Mvar. 

For the fault unspecified case in Section \ref{calECC}, instead of applying a specific fault, 
the perturbation in (\ref{pertb}) is applied to mimic the performance of dynamic var sources, for which 
$Q_1=Q_2$ and are chosen as 10 Mvar, 20 Mvar, 40 Mar, 80 Mvar, 160 Mvar, and 200 Mvar, $t_1=1$\,s, and $t_2=2$\,s.

\subsection{Case 1: Fault Specified}

\begin{figure}[!t]
\centering
\includegraphics[width=3.3in]{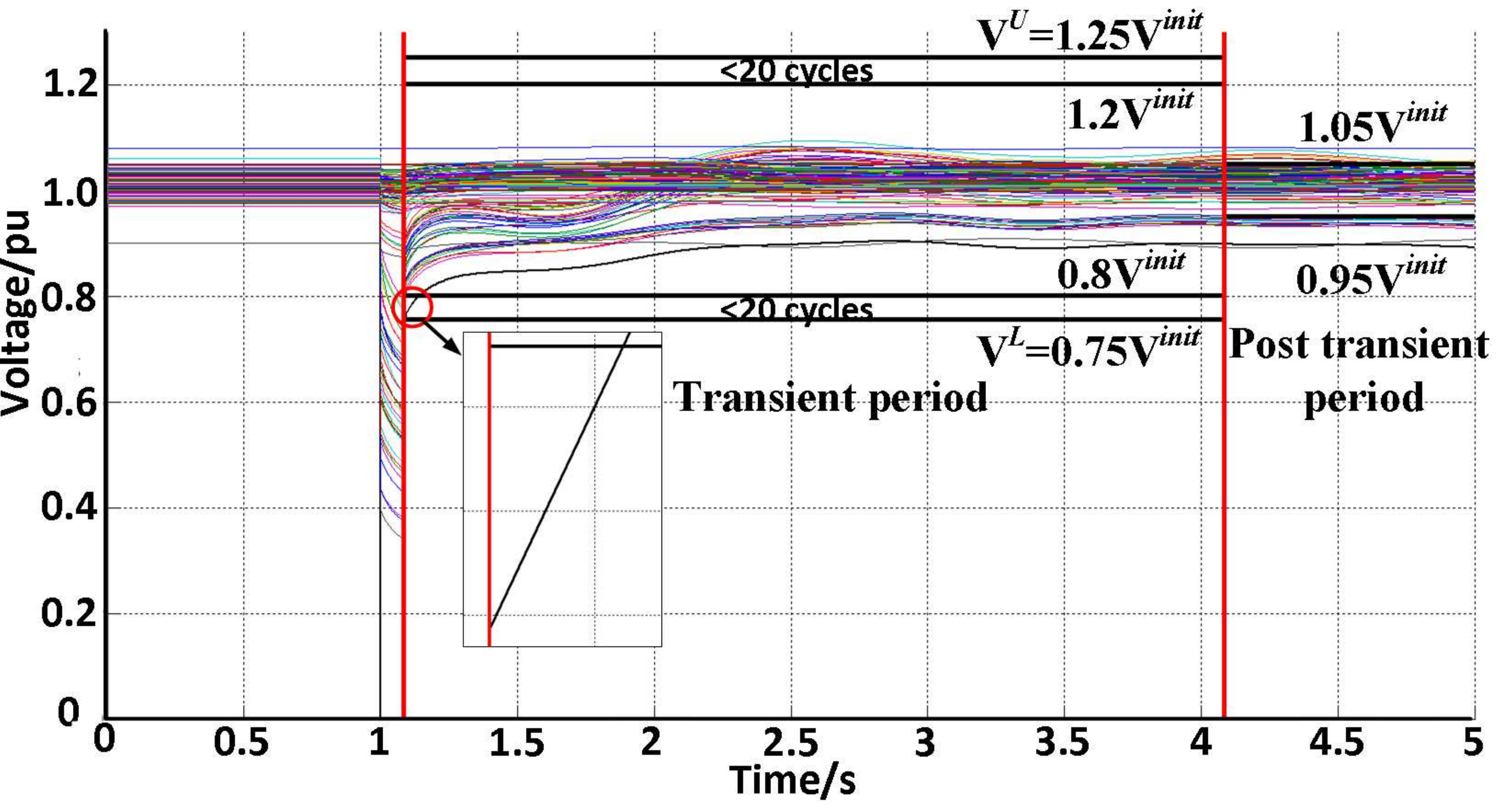}
\caption{Responses of all bus voltages under the most severe contingency.}
\label{FIDVR}
\end{figure}

\begin{figure}[!t]
\centering
\includegraphics[width=3.3in]{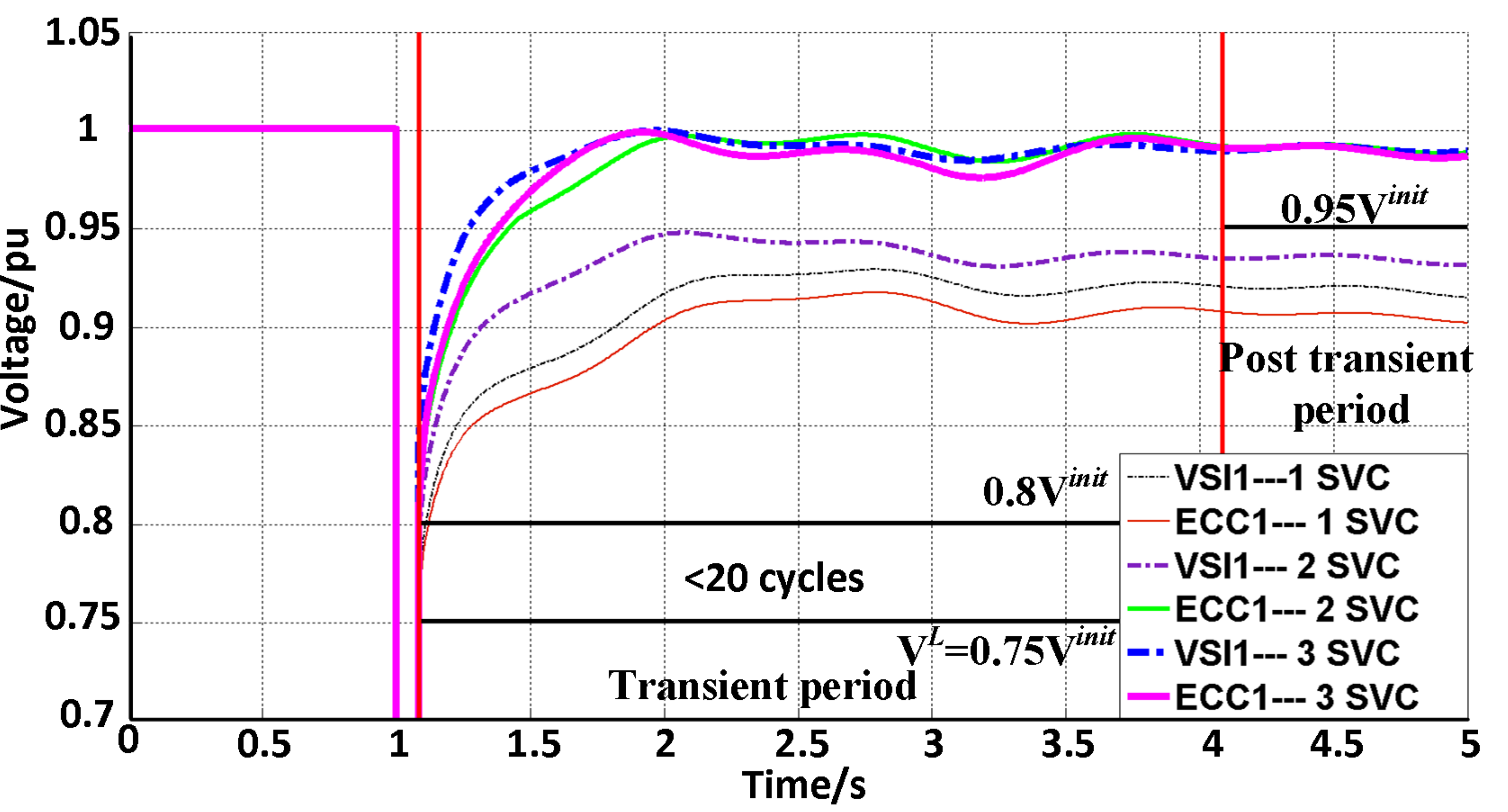}
\caption{Responses of most severe buses for placing 1 to 3 SVCs for both VSI1 and ECC1 methods.}
\label{FIDVR2}
\end{figure}

The most severe N-1 contingency is identified as a three-phase fault on bus 6 cleared by opening line $6-7$ after 5 cycles. 
The post-fault voltage trajectories are shown in Fig. \ref{FIDVR}. 
One bus violates the 25\% deviation limit at the beginning of the 3 seconds transient period 
and some buses violate the 5\% deviation limit at the post transient period.

First, the VSI-based method is applied to calculate the $\textrm{VSI}_i$ for each bus under the most severe contingency. 
It is called `VSI1' since only the most severe contingency is considered.
The top-20 buses in order are 3, 31, 6, 30, 34, 35, 36, 9, 16, 18, 17, 7, 12, 15, 14, 19, 20, 11, 114, and 110. 
If $v$ SVCs are installed, they should be placed at the top-$v$ buses in the list. 
Then, for the most severe contingency, the ECC-based method for which ECC is calculated by the method in Case 1 of Section \ref{calECC}, denoted by `ECC1', 
is applied to optimize the SVC locations. 
The locations for placing one to three SVCs by both methods are given in Table \ref{three}.

\begin{table}[!t]
\renewcommand{\arraystretch}{1.2}
\caption{Comparison of VSI1 Method and Fault Specified ECC1 Method}
\label{three}
\centering
\begin{tabular}{cp{0.4cm}p{0.4cm}p{0.6cm}p{0.6cm}p{1cm}p{1cm}}
\hline
\# of SVCs
& \multicolumn{2}{c}{One SVC}
& \multicolumn{2}{c}{Two SVCs}
&\multicolumn{2}{c}{Three SVCs}\\
\hline
 Method  & VSI1 & ECC1 & VSI1 & ECC1 & VSI1 & ECC1 \\
\hline
Optimal buses & 3 & 30 & 3, 31 & 6, 19 & 3, 31, 6 & 6, 12, 18 \\
FIDVR solved? & No  & No   & No  & Yes & Yes  & Yes  \\
\hline 
\end{tabular}
\end{table}

If only one SVC is installed, the placements from both methods cannot eliminate the FIDVR issue caused by the most severe contingency. 
If two SVCs are installed, the optimal placement from the ECC method solves the FIDVR issue while the VSI1 method does not. 
With three SVCs installed, the optimal placements by both methods can solve the issue. 
For the placements of one to three SVCs from both methods, the post-fault voltage responses of the most severe bus are shown in Fig. \ref{FIDVR2}. 
The violations happen on the lower limit 0.95 pu for the post-transient period.

\subsection{Case 2: Fault Unspecified}

Here, we respectively place 5 to 45 SVCs for comparison. 
The detailed placements of 5 to 40 SVCs for both VSI-based method and the fault unspecified ECC-based method (denoted by `ECC2') are listed in Table \ref{Placement}.

\begin{figure}[!t]
\centering
\includegraphics[width=3.3in]{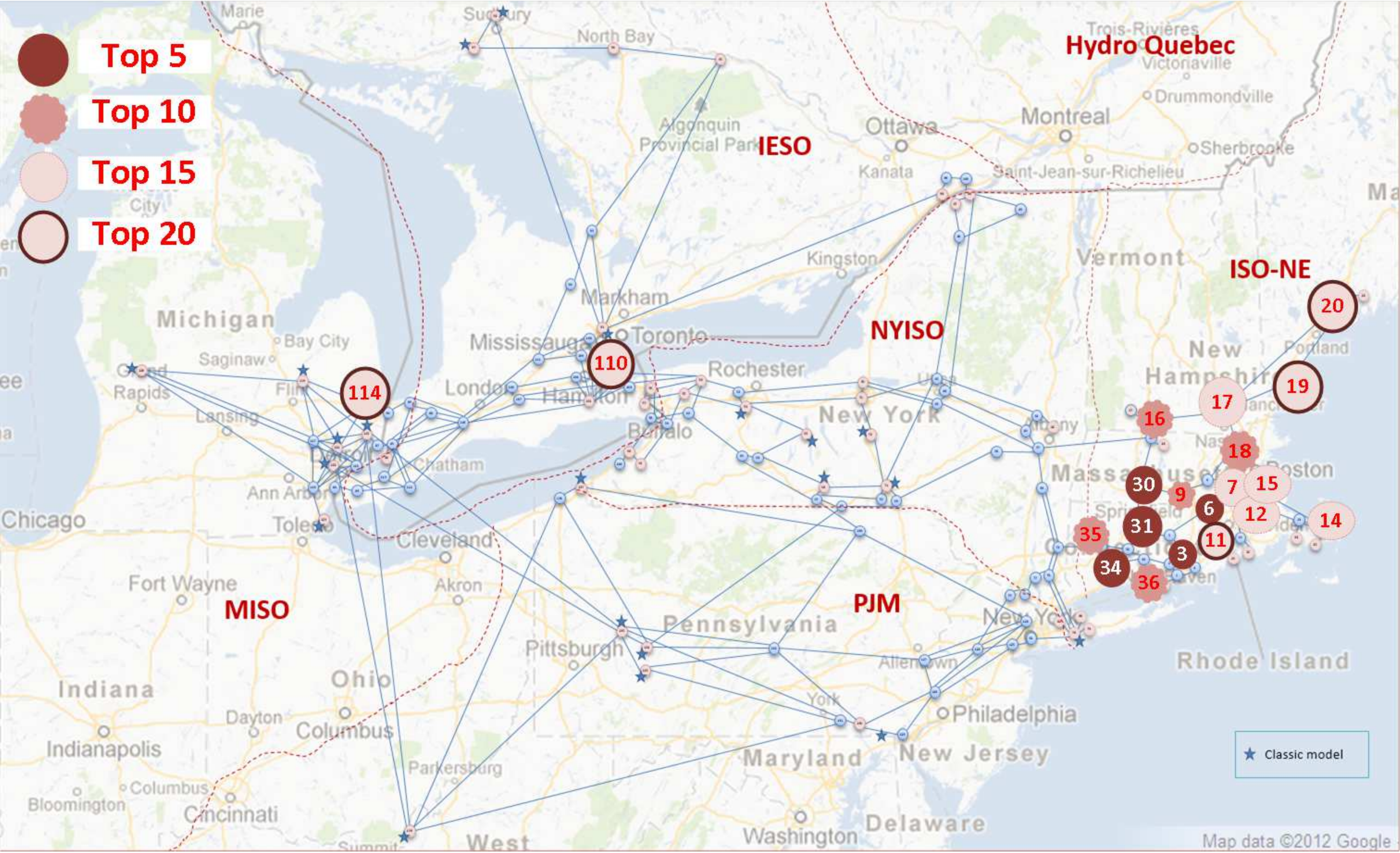}
\caption{Optimal 5, 10, 15, and 20 SVCs placed by VSI1 method.}
\label{NPCC-VSI1}
\end{figure}

\begin{figure}[!t]
\centering
\includegraphics[width=3.3in]{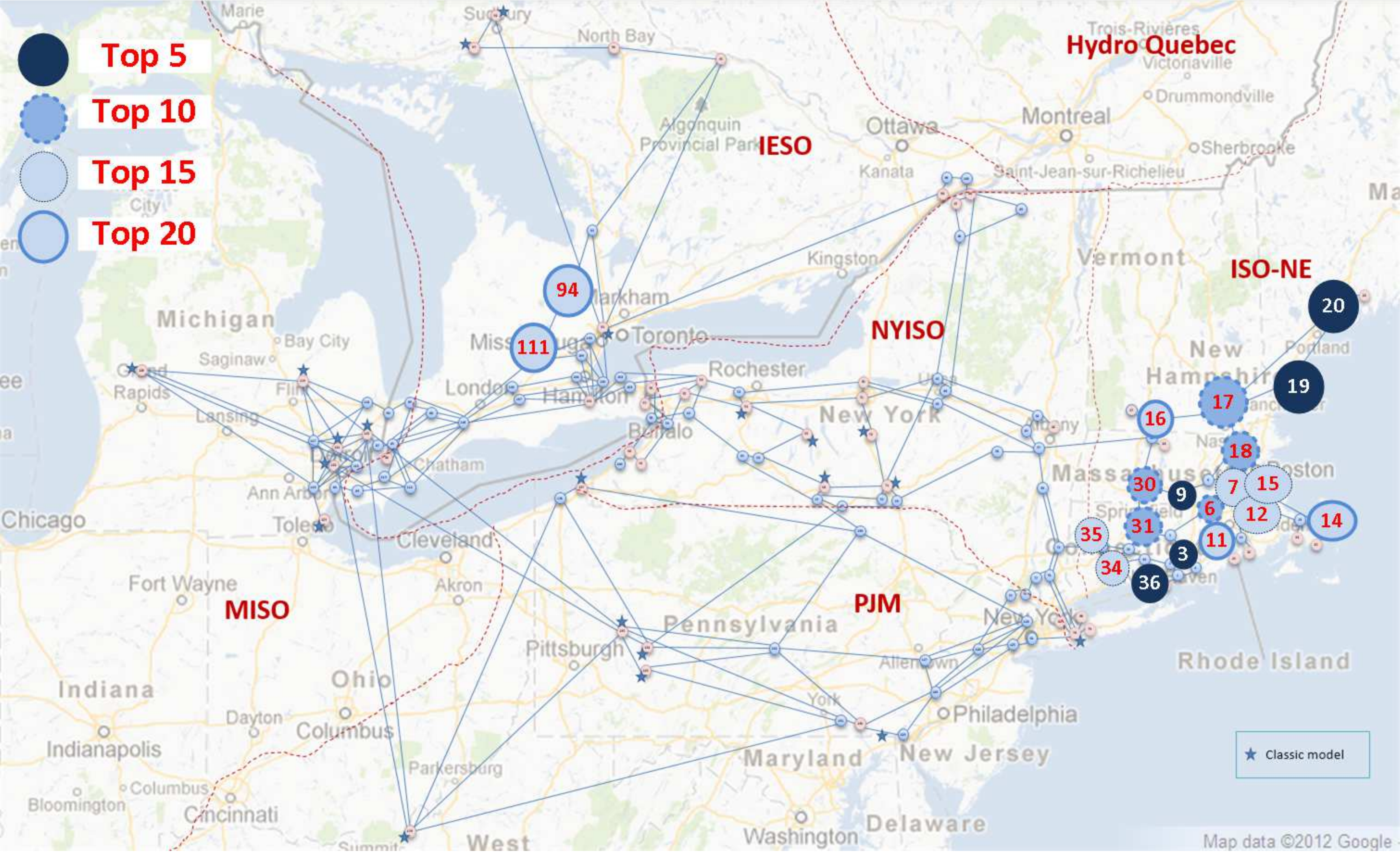}
\caption{Optimal 5, 10, 15, and 20 SVCs placed by VSI3 method.}
\label{NPCC-VSI3}
\end{figure}

The VSI-based method ranks potential locations based on one or several most severe contingencies.
The VSI-based method is called `VSI3' when the top three most severe contingencies are considered which include a three-phase fault at bus $6$ cleared by opening line $6$--$7$, 
a three-phase fault at bus $17$ cleared by opening line $17$--$16$, and a three-phase fault at bus $16$ cleared by opening line $16$--$29$.
Figs. \ref{NPCC-VSI1}--\ref{NPCC-VSI3} illustrates the results for placing 5, 10, 15, and 20 SVCs by the VSI1 and VSI3 method. 

\begin{figure}[!t]
\centering
\includegraphics[width=3.3in]{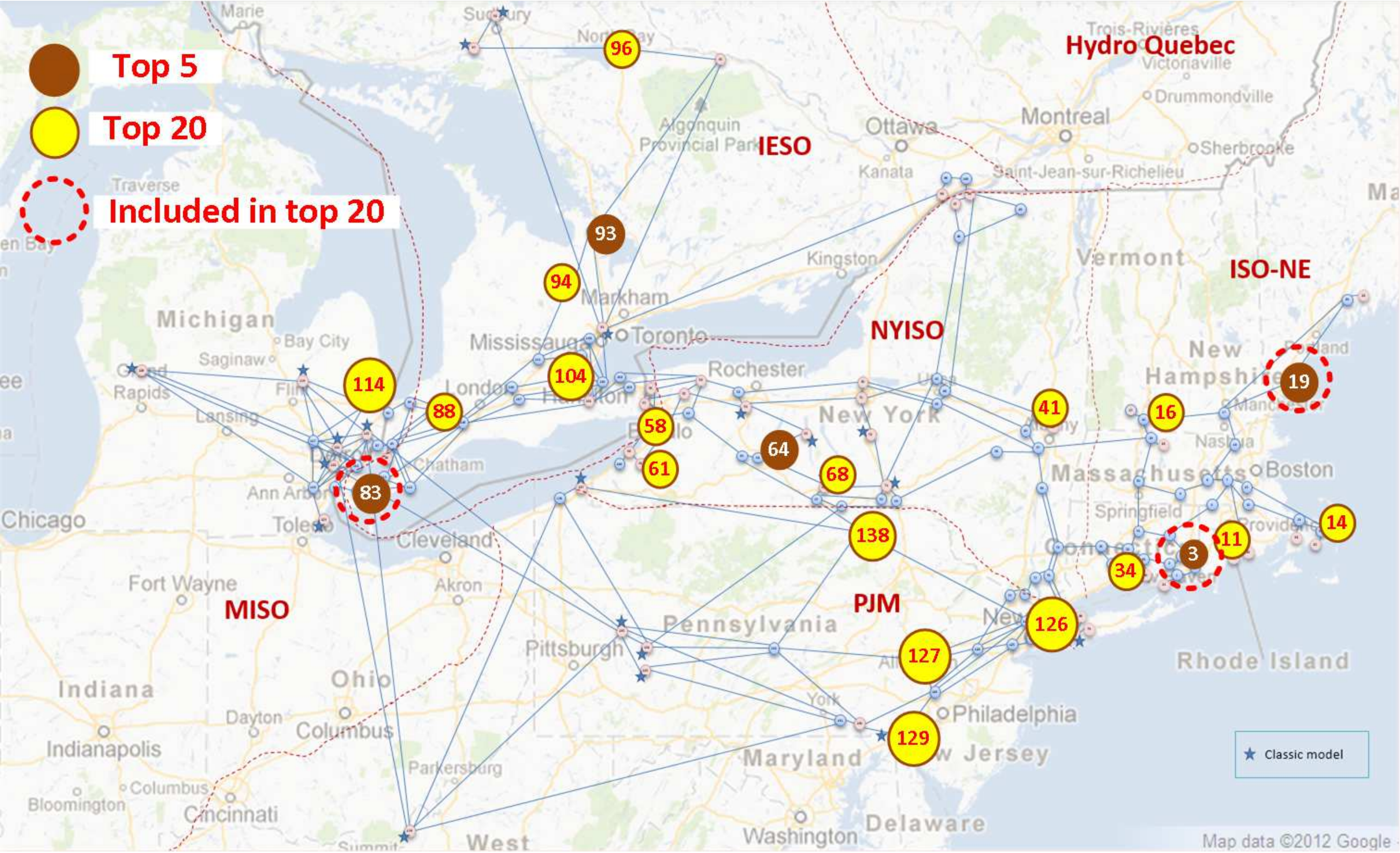}
\caption{Optimal 5 and 20 SVCs placed by the fault unspecified ECC2 method.}
\label{NPCC-ECC}
\end{figure}

However, that is not the case for the ECC-based method, which optimizes the placement for each number of SVCs individually without knowing the specific fault information. 
Fig. \ref{NPCC-ECC} illustrates the optimal placements for 5 and 20 SVCs by the ECC2 method, for which only buses 3, 19, and 83 appear in both top 5 and top 20 placements.

\begin{table}[!t]
\renewcommand{\arraystretch}{1.0}
\caption{SVC Placements from VSI-Based Method and Fault Unspecified ECC-Based Method}
\label{Placement}
\centering
\begin{tabu}{ccc}
\hline
\# of SVCs & Method & Optimal placement \\
\hline
\multirow{2}{*}{5} & VSI1 & 3, 31, 6, 30, 34  \\
  & VSI3 & 3, 20, 19, 36, 9 \\
  & ECC2 & 3, 19, 64, 83, 93 \\
\hline
\multirow{2}{*}{10} & VSI1 & 3, 31, 6, 30, 34, 35, 36, 9, 16, 18  \\
  & VSI3 & 3  20, 19, 36, 9, 6, 30, 17, 31, 18 \\
  & ECC2 & 3, 14, 19, 46, 61, 83, 93, 95, 124, 138 \\
\hline 
\multirow{3}{*}{15} & VSI1 & \tabincell{l}{3, 31, 6, 30, 34, 35, 36, 9, \\ 16, 18, 17, 7, 12, 15, 14} \\
  & VSI3 & \tabincell{l}{3, 20, 19, 36, 9, 6, 30, 17, \\ 31, 18, 34, 35, 15, 7, 12} \\
  & ECC2 & \tabincell{l}{3, 11, 14, 19, 34, 61, 64, 83, \\ 89, 93, 96, 113, 126, 129, 138}  \\
\hline
\multirow{3}{*}{20} & VSI1 & \tabincell{l}{3, 31, 6, 30, 34, 35, 36, 9, 16, 18, 17, \\ 7, 12, 15, 14, 19, 20, 11, 114, 110}  \\ 
  & VSI3 & \tabincell{l}{3, 20, 19, 36, 9, 6, 30, 17, 31, 18, 34, \\ 35, 15, 7, 12, 16, 14, 11, 111, 94} \\
  & ECC2 & \tabincell{l}{3, 11, 14, 16, 19, 34, 41, 58, 61, 68, 83, \\ 88, 94, 96, 104, 114, 126, 127, 129, 138} \\
\hline
\multirow{4}{*}{25} & VSI1 & \tabincell{l}{3, 31, 6, 30, 34, 35, 36, 9, 16, 18, \\ 17, 7, 12, 15, 14, 19, 20, 11, 114,\\ 110, 90, 104, 89, 109, 108}  \\
  & VSI3 & \tabincell{l}{3, 20, 19, 36, 9, 6, 30, 17, 31, 18, \\ 34, 35, 15, 7, 12, 16, 14, 11, 111, \\ 94, 113, 93, 90, 83, 89} \\
  & ECC2 & \tabincell{l}{3, 9, 11, 12, 14, 16, 19, 34, 41, 45, \\ 58, 61, 68, 83, 88, 94, 95, 104, 113, \\ 114, 124, 126, 127, 129, 138}  \\
\hline
\multirow{4}{*}{30} & VSI1 & \tabincell{l}{3, 6, 7, 9, 11, 12, 14, 15, 16, 17, 18, 19, \\ 20, 30, 31, 34, 35, 36, 41, 56, 89, 90, \\ 104, 106, 107, 108, 109, 110, 111, 114}  \\
  & VSI3 & \tabincell{l}{3, 20, 19, 36, 9, 6, 30, 17, 31, 18, 34, 35, \\ 15, 7, 12, 16, 14, 11, 111, 94, 113, 93, \\ 90, 83, 89, 114, 106, 109, 108, 105} \\
  & ECC2 & \tabincell{l}{3, 6, 9, 11, 12, 14, 18, 19, 34, 41, 46, 59, \\ 61, 64, 80, 83, 88, 93, 94, 95, 96, 104, \\ 113, 114, 124, 126, 127, 129, 136, 138}  \\
\hline
\multirow{5}{*}{35} & VSI1 & \tabincell{l}{3, 6, 7, 9, 11, 12, 14, 15, 16, 17, 18, \\ 19, 20, 30, 31, 34, 35, 36, 41, 42, 45, \\ 56, 74, 89, 90, 93, 94, 104, 106, 107, \\ 108, 109, 110, 111, 114}  \\
  & VSI3 & \tabincell{l}{3, 20, 19, 36, 9, 6, 30, 17, 31, 18, 34, \\ 35, 15, 7, 12, 16, 14, 11, 111, 94, 113, \\ 93, 90, 83, 89, 114, 106, 109, 108, \\ 105, 107, 104, 88, 110, 41} \\
  & ECC2 & \tabincell{l}{3, 6, 9, 11, 12, 14, 15, 18, 19, 30, 34, \\ 41, 46, 53, 59, 61, 64, 80, 83, 88, 93, \\ 94, 95, 96, 104, 113, 114, 119, 124, \\ 125, 126, 127, 129, 136, 138}  \\
\hline
\multirow{5}{*}{40} & VSI1 & \tabincell{l}{3, 6, 7, 9, 11, 12, 14, 15, 16, 17, 18, 19, \\ 20, 30, 31, 34, 35, 36, 41, 42, 45, 56, 74, \\ 78, 86, 89, 90, 93, 94, 95, 104, 106, 107, \\ 108, 109, 110, 111, 114, 131, 132}  \\
  & VSI3 & \tabincell{l}{3, 20, 19, 36, 9, 6, 30, 17, 31, 18, 34, 35, \\ 15, 7, 12, 16, 14, 11, 111, 94, 113, 93, 90, \\ 83, 89, 114, 106, 109, 108, 105, 107, 104, \\ 88, 110, 41, 42, 95, 74, 116, 45} \\
  & ECC2 & \tabincell{l}{3, 6, 9, 11, 12, 14, 15, 18, 19, 30, 31, 34, \\ 41, 46, 53, 59, 61, 64, 68, 80, 83, 88, 89, \\ 90, 93, 94, 95, 96, 104, 113, 114, 119, \\ 124, 125, 126, 127, 128, 129, 136, 138}  \\
\hline  
\end{tabu}
\end{table}

The VSI-based method depends on the specific setting of the contingency, such as the fault duration, while the fault unspecified ECC2 method does not. 
Thus it is expected that the ECC2 method has better performance under different contingency settings. 
In Table \ref{FIDVR-Coverage}, the VSI1, VSI3, and ECC2 methods are compared in terms of how many contingencies with FIDVR issues are addressed with 4, 5, and 6 cycles of fault duration. 
As mentioned above, for the NPCC system 40 contingencies with a fault duration of 4 to 6 cycles can cause FIDVR issues. 
From Table \ref{FIDVR-Coverage}, it is seen that the ECC2 method can address more FIDVR issues in general than the VSI-based method. 
For example, for 5-cycle fault duration, the ECC2 method can address 28 out of the 40 contingencies (70\% coverage)   
by placing 40 SVCs while the VSI1 and VSI3 methods can only address 25 (62.5\% coverage) and 24 contingencies (60\% coverage).

\begin{table}[!t]
\renewcommand{\arraystretch}{1.0}
\caption{Percentage Coverage for Different Fault Durations}
\label{FIDVR-Coverage}
\centering
\begin{tabular}{ccccc}
\hline
\multirow{2}{*}{Cycle} & \multirow{2}{*}{\# of SVCs} & \multicolumn{3}{c}{Total FIDVR contingencies addressed} \\
               &              &   VSI1 & VSI3 & ECC2 \\
\hline
		 & 5 & 14 & 12 & 15 \\
         & 10 & 14 & 14 & 18 \\
         & 15 & 16 & 16 & 22 \\
         & 20 & 21 & 18 & 27 \\
4        & 25 & 21 & 21 & 28 \\
         & 30 & 23 & 21 & 28 \\
         & 35 & 25 & 28 & 28 \\
         & 40 & 26 & 28 & 28 \\
         & 45 & 26 (65\%) & 28 (70\%) & 28 (70\%)\\
\hline
		 & 5 & 10 & 10 & 10 \\
         & 10 & 13 & 13 & 16 \\
         & 15 & 15 & 15 & 20 \\
         & 20 & 19 & 18 & 26 \\
5        & 25 & 19 & 20 & 27 \\
         & 30 & 22 & 20 & 27 \\
         & 35 & 22 & 23 & 27 \\
         & 40 & 25 & 24 & 28 \\
         & 45 & 25 (62.5\%) & 24 (60\%) & 28 (70\%)\\
\hline 
		 & 5 & 3 & 4 & 6\\
         & 10 & 9 & 11 & 10\\
         & 15 & 15 & 15 & 16\\
         & 20 & 19 & 17 & 21\\
6        & 25 & 19 & 20 & 24 \\
         & 30 & 21 & 20 & 24\\
         & 35 & 24 & 22 & 26\\
         & 40 & 25 & 24 & 26\\
         & 45 & 25 (62.5\%) & 24 (60\%) & 27 (67.5\%) \\
\hline
\end{tabular}
\end{table}

\begin{figure}[!t]
\centering
\includegraphics[width=2.5in]{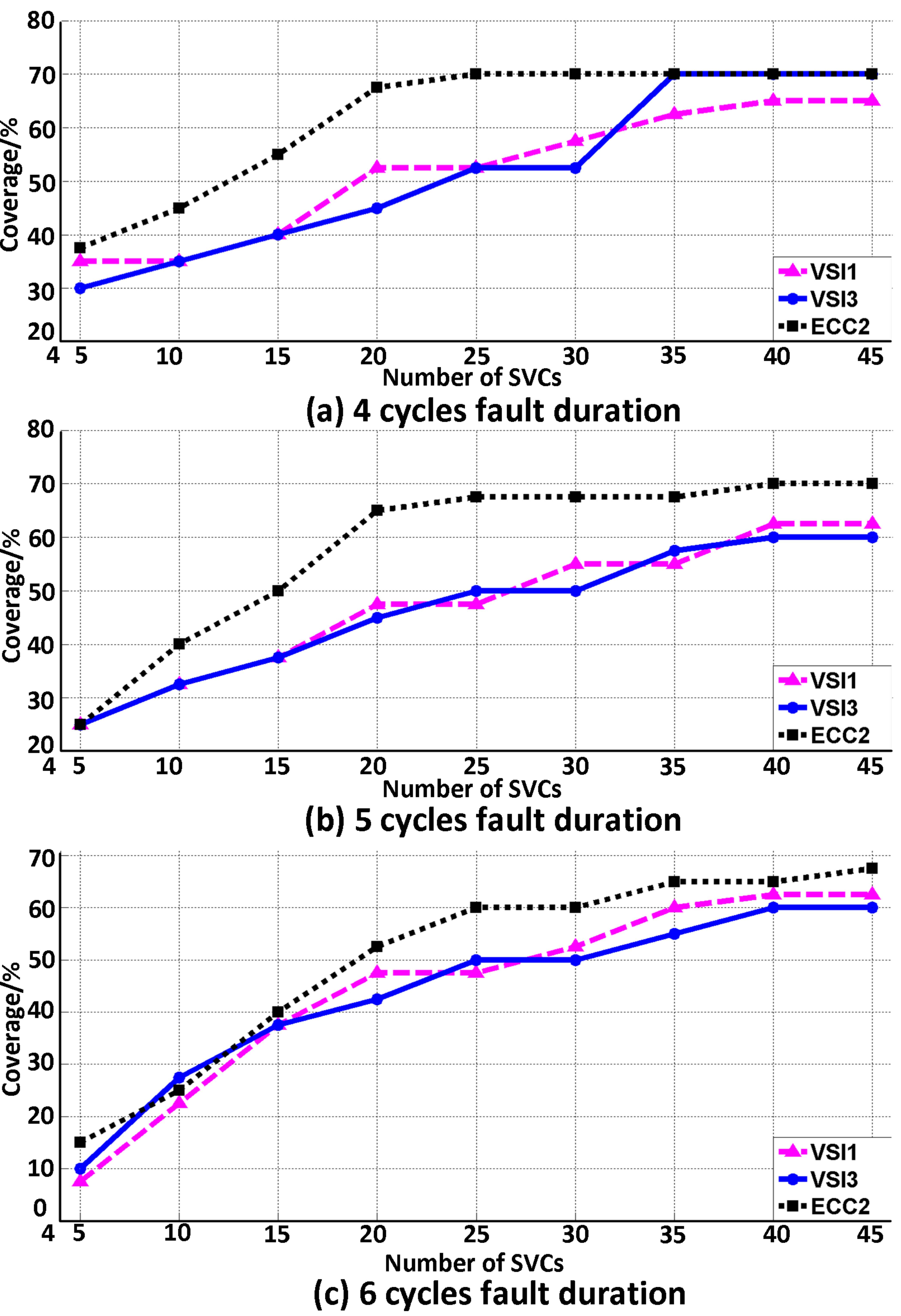}
\caption{Percentage coverages for VSI-based method and fault unspecified ECC2 method under different fault durations.}
\label{Trend}
\end{figure}

\begin{figure}[!t]
\centering
\includegraphics[width=2.7in]{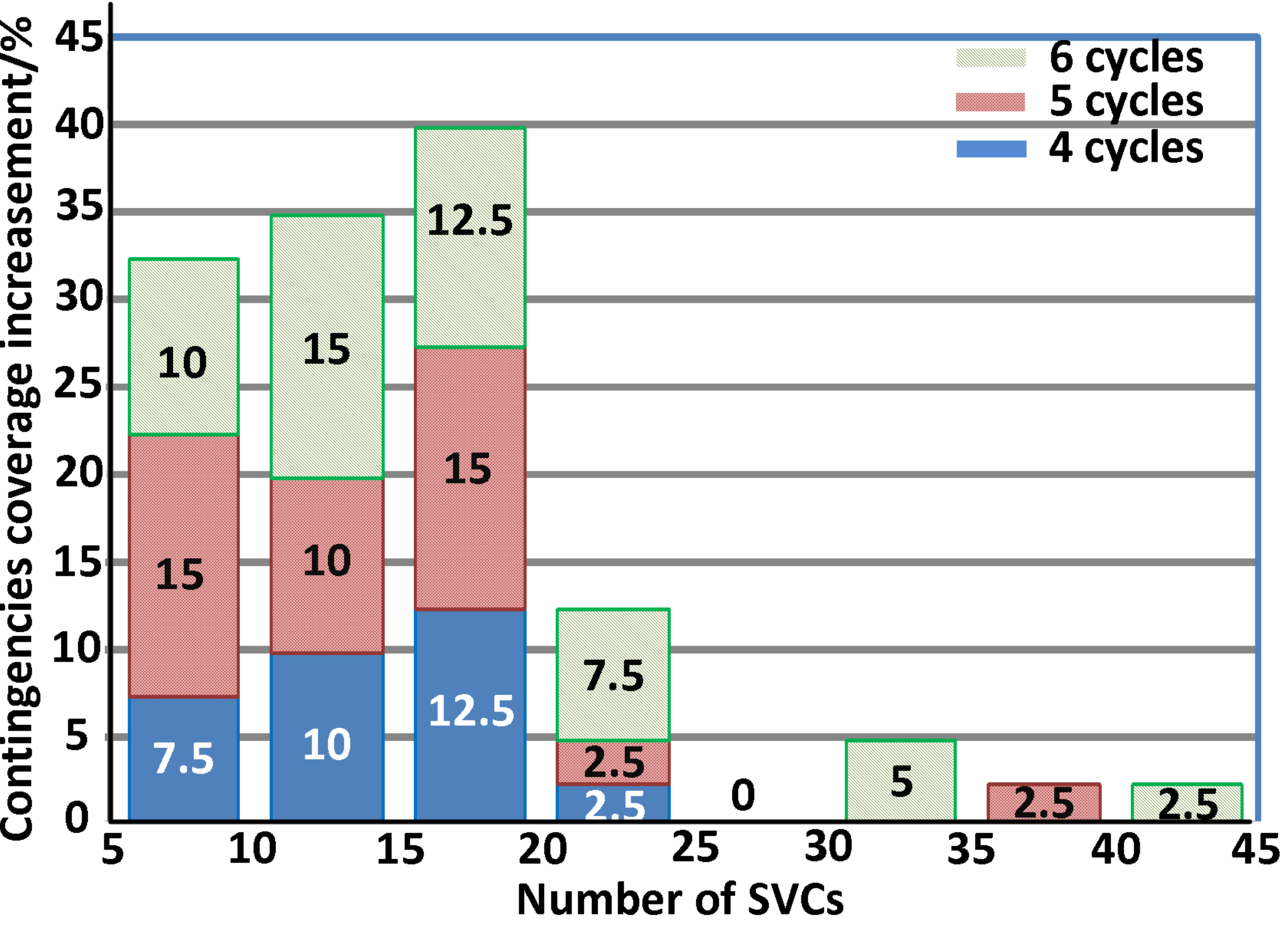}
\caption{Percentage coverage increment for fault unspecified ECC2 method under different fault durations.}
\label{fig11}
\end{figure}

Fig. \ref{Trend} shows the percentage coverage of the contingencies with FIDVR issues addressed by the VSI-based and fault unspecified ECC2 methods. 
Fig. \ref{fig11} shows the increase of the percentage coverage with the increase of the number of SVCs for the  ECC2 method. 
For example, from 5 to 10 SVCs installation, the coverage increases by 7.5\%, 15\%, and 10\% for 4-, 5-, and 6-cycle fault duration, respectively. 
It is seen that the coverage grows rapidly when the number of SVCs is less than 25, and then does not grow for 4-cycle fault duration and hardly grow for 5- and 6-cycle duration. 
From Table \ref{FIDVR-Coverage} and Figs. \ref{Trend} and \ref{fig11}, it is also seen that there is a maximum number of SVCs above which the contingency coverage will not increase. 
For example, for 5-cycle fault duration, this maximum number is 40.

The optimal number of SVCs that minimizes the total cost can also be determined.  
The total cost includes that for installing SVCs ($C_{\textrm{SVC}}$) and that for the FIDVR issues caused by the contingencies that are not addressed for the considered fault durations ($C_{\textrm{FIDVR}}$). 
The total cost can be written as:   
\begin{align}
C&=C_{\textrm{SVC}} + C_{\textrm{FIDVR}} \notag \\
 &=c_{\textrm{SVC}}\,n_{\textrm{SVC}} + c_{\textrm{FIDVR}} \sum\limits_{i=1}^{K} \big(N_{\textrm{cont}} - n^i_{\textrm{cont}}(n_{\textrm{SVC}})\big)
\end{align}
where $c_{\textrm{SVC}}$ is the cost of one SVC, $n_{\textrm{SVC}}$ is the number of SVCs installed, 
$c_{\textrm{FIDVR}}$ is the cost of FIVDR issue caused by one contingency, $N_{\textrm{cont}}=40$ is the total number of 
contingencies that are considered, and $n^i_{\textrm{cont}}$ is the number of contingencies addressed by installing 
$n_{\textrm{SVC}}$ SVCs which is a function of $n_{\textrm{SVC}}$, and $K=3$ is the number of fault durations. 

With the increase of $n_{\textrm{SVC}}$, $C_{\textrm{SVC}}$ will linearly increase while $C_{\textrm{FIDVR}}$ will decrease before the effect of installing SVCs in eliminating the contingencies with FIDVR issues saturates. Therefore, we can get an optimal number of SVCs ($n^*_{\textrm{SVC}}$) that minimizes the total cost $C$. 
By assuming $c_{\textrm{FIDVR}}=5\,c_{\textrm{SVC}}$, the optimal number of SVCs is determined as 25 and the corresponding minimal cost is $230 \, c_{\textrm{SVC}}$.

\subsection{Voltage Collapse Scenario}

\begin{figure}[!t]
\centering
\includegraphics[width=2.3in]{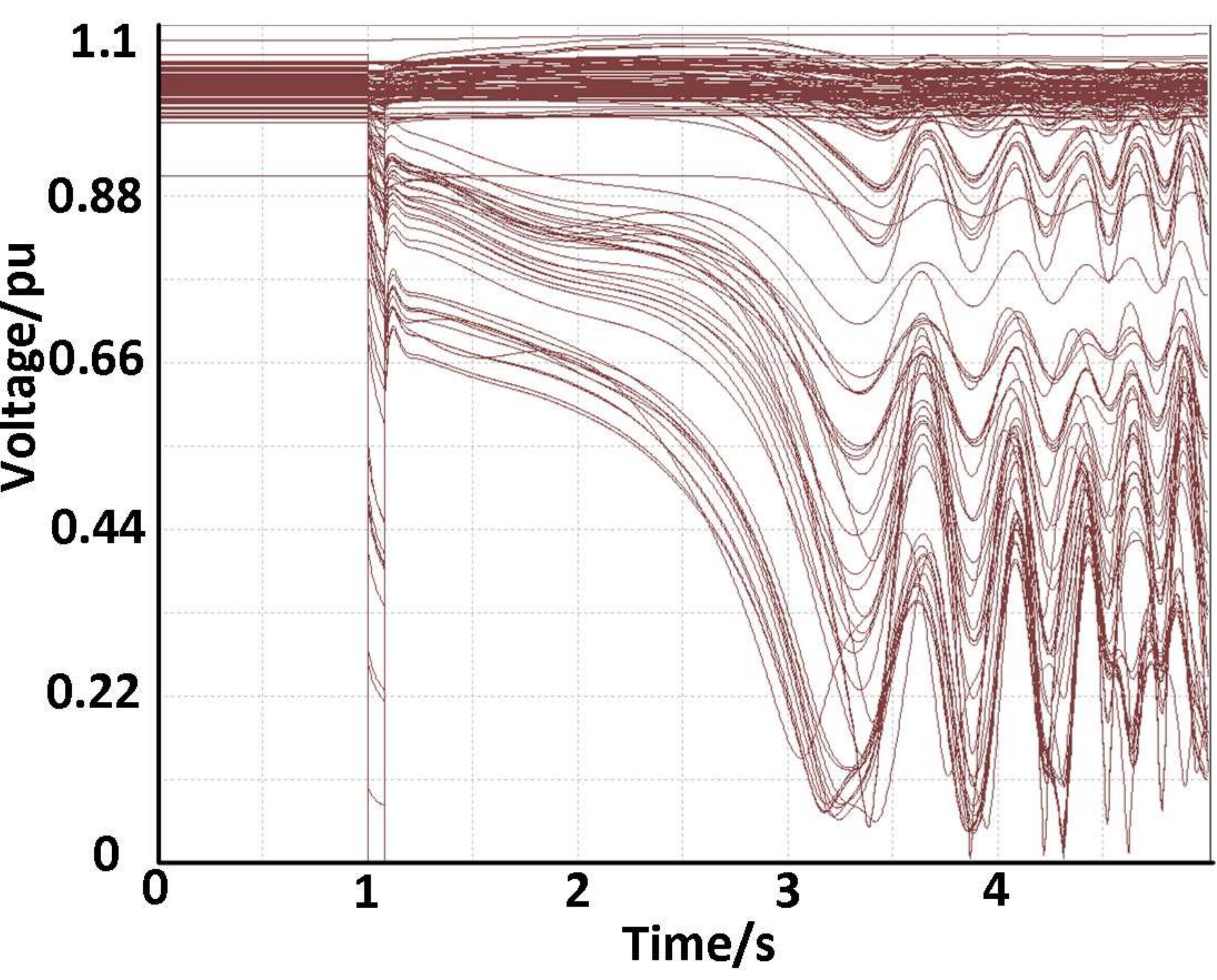}
\caption{Voltage recovery without var support under voltage collapse case.}
\label{fig12}
\end{figure}

\begin{figure}[!t]
\centering
\includegraphics[width=3.3in]{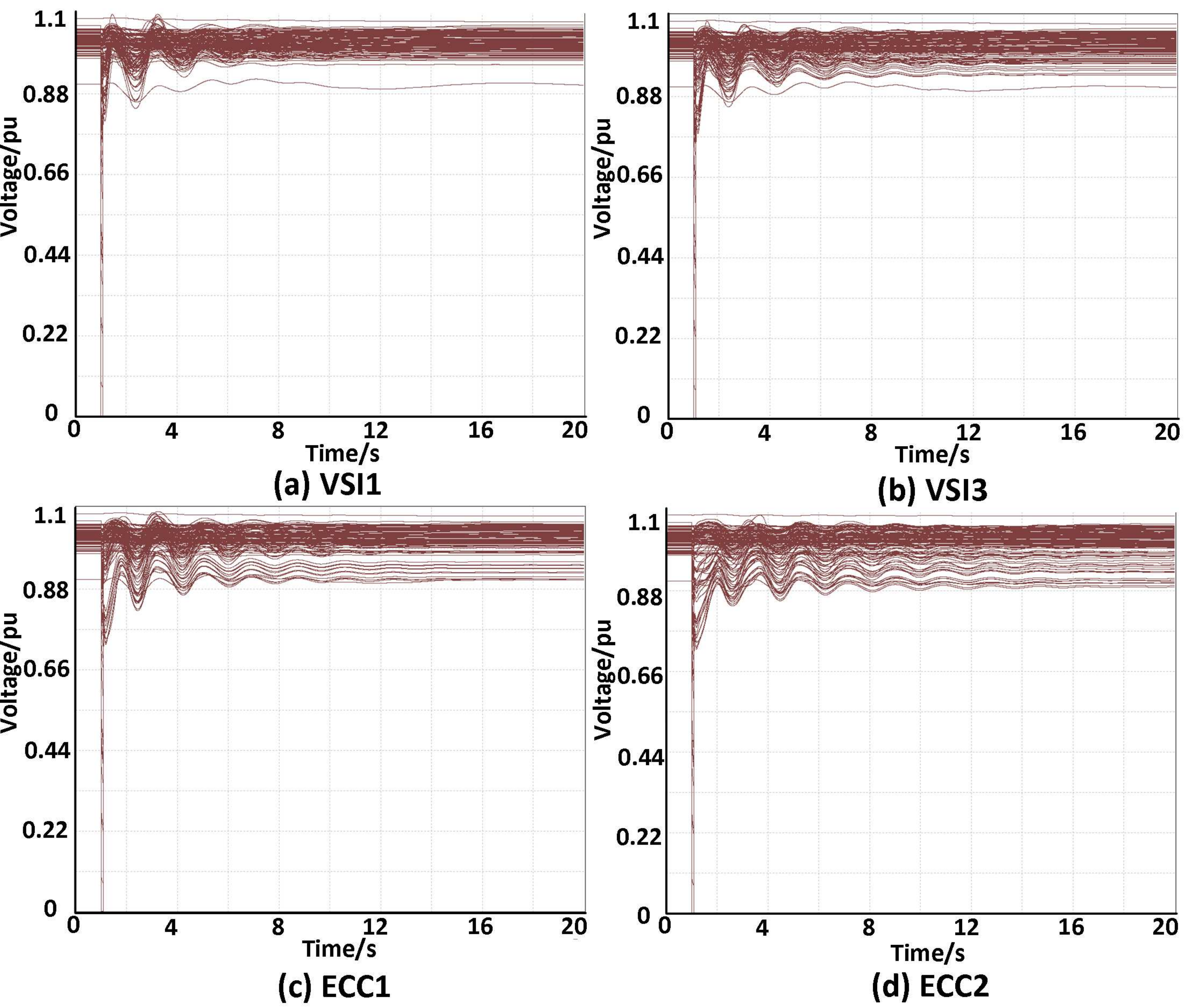}
\caption{Voltage recovery for four methods installing 5 SVCs.}
\label{fig13}
\end{figure}

We also evaluate the VSI1, VSI3, ECC1, and ECC2 methods under an extreme case, in which a three-phase fault is applied at bus 35 and is cleared by opening line $35-73$ after 5 cycles, 
and a half of the generation at bus 21 is lost when the line is tripped. 
As shown in Fig. \ref{fig12}, the voltages collapse after the fault is cleared. 
Fig. \ref{fig13} illustrates the voltage recovery under the four methods with 5 SVCs installed. 
The SVC placement for ECC1 is bus 6, 12, 19, 30, and 56 while those for VSI1, VSI3, and ECC2 can be found in Table \ref{Placement}.

From Fig. \ref{fig13} it is seen that the system voltages can be successfully maintained against this voltage collapse scenario by 5 SVCs placed by all four methods. 
Since the location of the initiating fault of the voltage collapse is close to the most severe contingency used by ECC1 and VSI1 and also the top three most severe contingencies used by VSI3, 
the SVCs placed by them are closer to the initiating fault and are more effective in preventing voltage collapse and maintaining system stability.
Because ECC2 is fault independent and focuses on the controllability of voltages across the whole system, it is not as effective as the other methods. 
The resulting damping ratios of the system with 5 SVCs placed by VSI1, VSI3, ECC1, and ECC2 are estimated by the Prony analysis as 10.20\%, 8.07\%, 9.51\%, and 6.67\%, respectively, 
which also verifies the above comparison of the four methods.

\subsection{Influence of the Capacity of SVCs}

\begin{figure}[!t]
\centering
\includegraphics[width=2.9in]{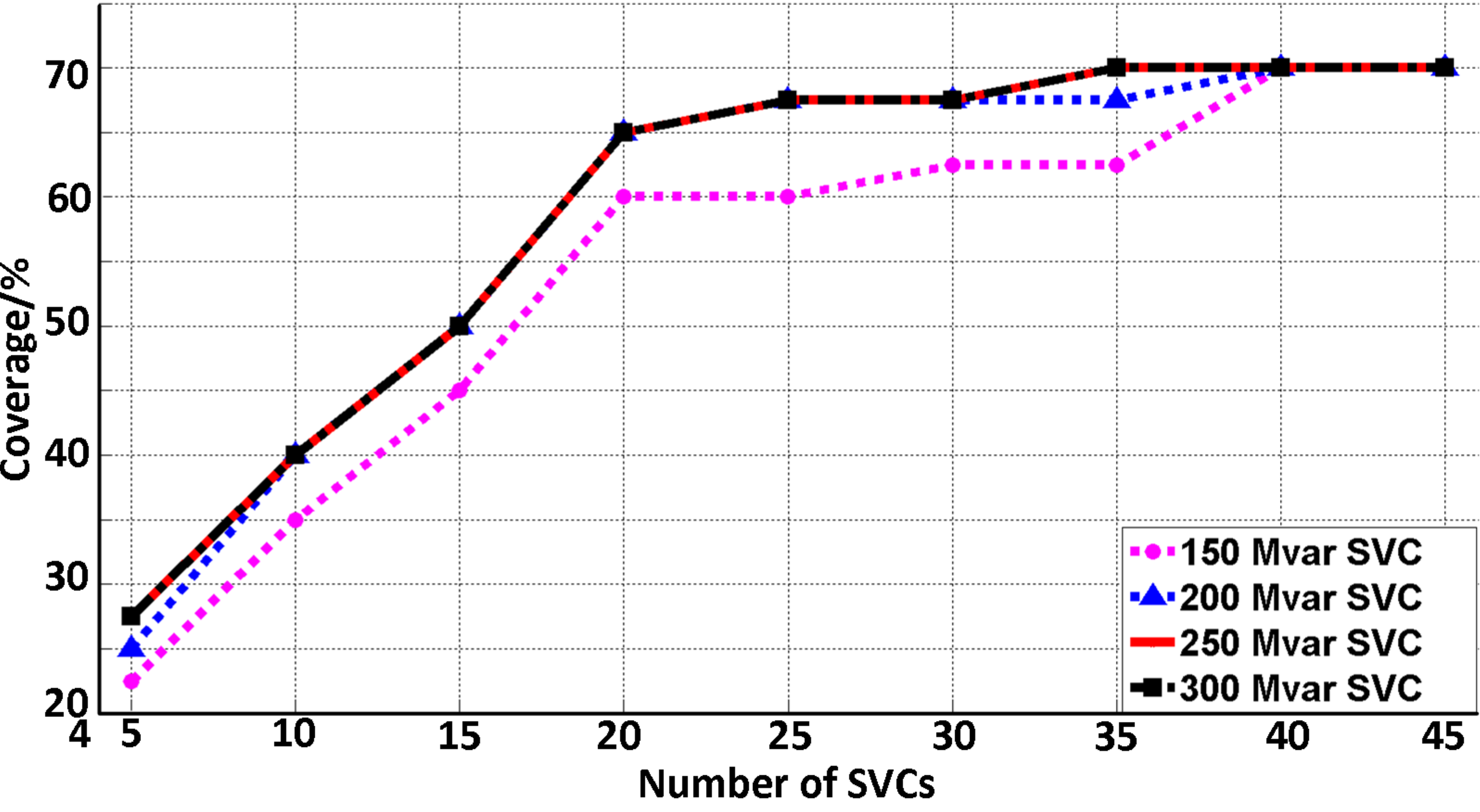}
\caption{Percentage coverages for ECC2 with different SVC capacities.}
\label{fig14}
\end{figure}

In the above sections we assume that all of the SVCs have an identical capacity of 200 Mvar. 
Here, we discuss how SVC sizes may influence the results. 
Fig. \ref{fig14} illustrates the percentage coverage for different SVC capacities. 
It is seen that the improvement of contingency coverages is not obvious when the SVC capacity is above 200 Mvar, 
which also implies that it is practical to install SVCs with a 200 Mvar capacity.

\section{Conclusion} \label{conclusion}

In this paper, ECC is applied to quantify the degree of controllability of the voltage magnitudes 
under a specific dynamic var source placement and the optimal placement problem for dynamic var sources is formulated as an optimization problem 
that maximizes the determinant of ECC.
 
The proposed method is tested and validated on an NPCC 140-bus system. 
The results show that the proposed method for fault specified case can solve the FIDVR issue caused by the most severe N-1 contingency 
with fewer dynamic var sources than that of the VSI-based method. 
The proposed method for fault unspecified case does not depend on the settings of the contingency and thus has 
better performance under different fault durations, in the sense that when placing the same number of SVCs 
the proposed method can address more FIDVR issues than the VSI-based method. 
The optimal number of SVCs that minimizes the total cost is determined and it is shown that the proposed method can help mitigate voltage collapse. 
It is also found that the improvement is not obvious when the SVC capacity is above 200 Mvar.

\begin{IEEEbiography} [{\includegraphics[width=1in,height=1.25in,clip,keepaspectratio]{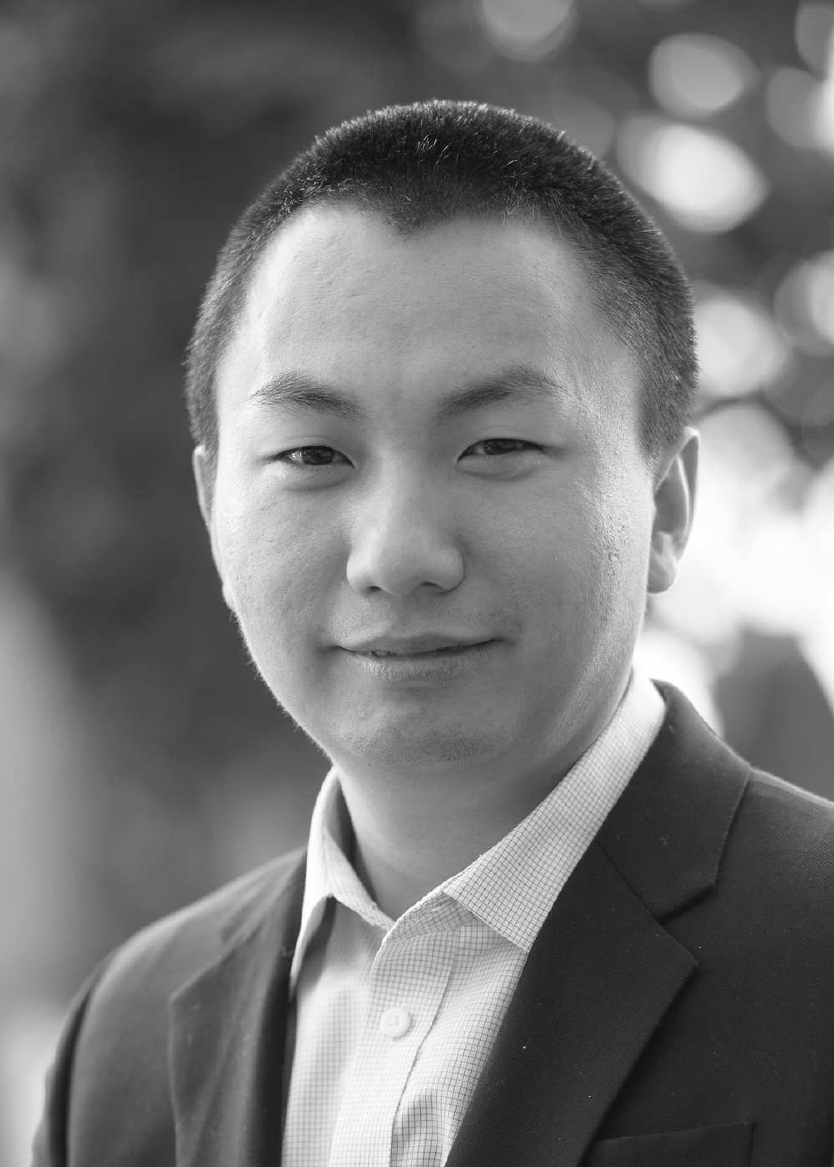}\vfill}]
{Junjian Qi} (S'12--M'13)
received the B.E. degree from Shandong University, Jinan, China, in 2008 and the Ph.D. degree Tsinghua University, Beijing, China, in 2013, both in electrical engineering.

In February--August 2012 he was a Visiting Scholar at Iowa State University, Ames, IA, USA. During September 2013--January 2015 he was 
a Research Associate at Department of Electrical Engineering and Computer Science, University of Tennessee, Knoxville, TN, USA. 
Currently he is a Postdoctoral Appointee at the Energy Systems Division, Argonne National Laboratory, Argonne, IL, USA. 
His research interests include cascading blackouts, power system dynamics, state estimation, synchrophasors, and cybersecurity.
\end{IEEEbiography}

\begin{IEEEbiography} [{\includegraphics[width=1in,height=1.25in,clip,keepaspectratio]{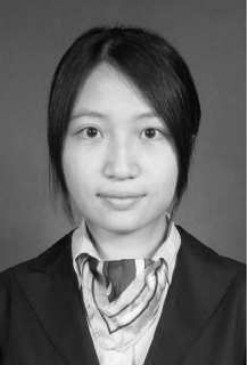}\vfill}]
{Weihong Huang} (S'14)
received the B.E degree from Huazhong University of Science and Technology, Wuhan, China, in 2009, and the master degree from University of Tennessee, Knoxville, TN, USA in 2014, 
both in electrical engineering.

She worked at Yueqing Electric Power Bureau, Grid State Corporation of China, Wenzhou, China, from 2009 to 2012. 
Currently, she is pursuing Ph.D. degree at Department of Electrical Engineering and Computer Science, University of Tennessee, Knoxville, TN, USA. 
Her research interests include operation control of dynamics reactive power sources.
\end{IEEEbiography}

\begin{IEEEbiography} [{\includegraphics[width=1in,height=1.25in,clip,keepaspectratio]{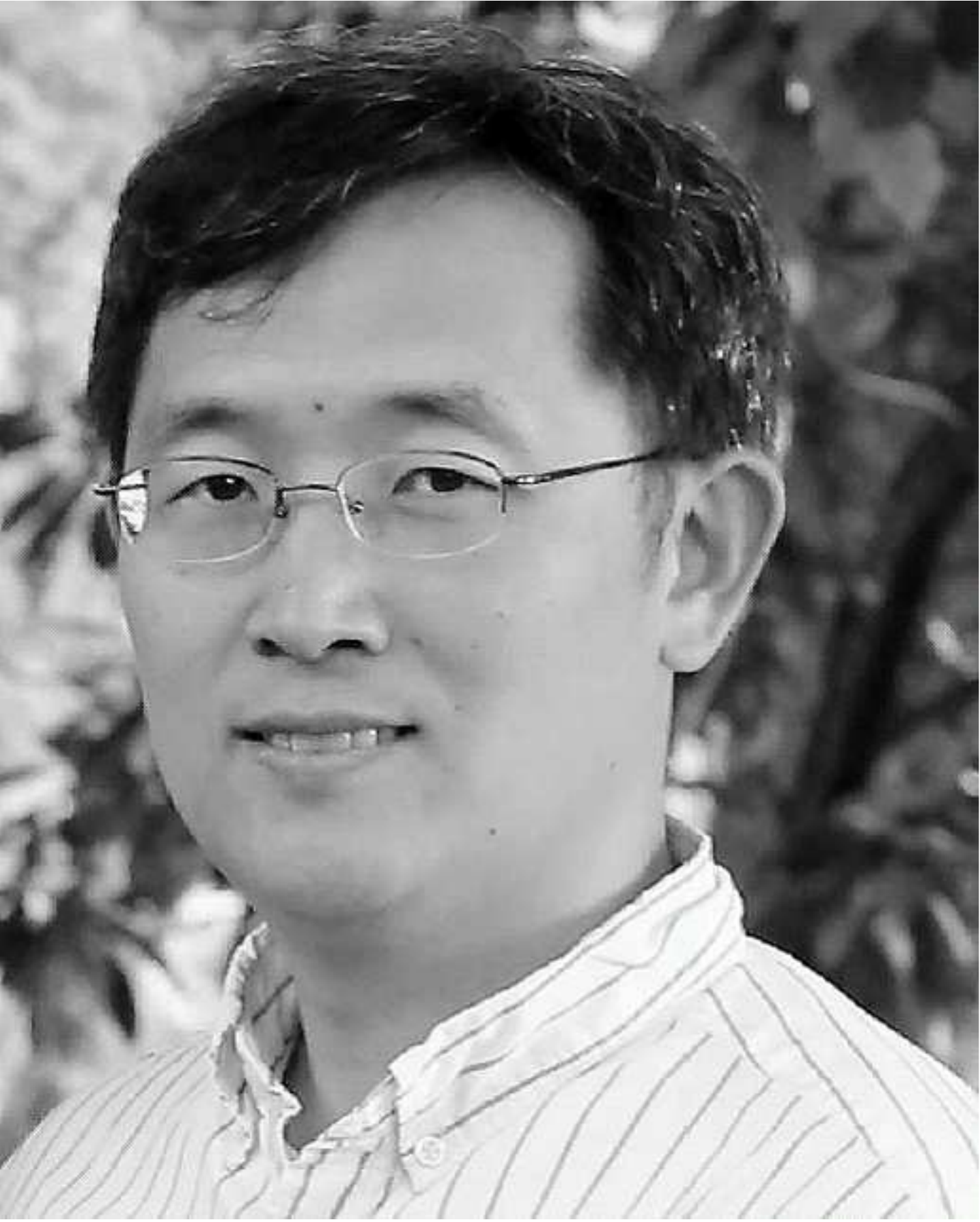}\vfill}]
{Kai Sun} (M'06--SM'13)
received the B.S. degree in automation in 1999 and the Ph.D. degree in control science and engineering in 2004 both from Tsinghua University, Beijing, China. 

He is currently an assistant professor at the Department of EECS, University of Tennessee in Knoxville. 
He was a project man-ager in grid operations and planning at the EPRI, Palo Alto, CA from 2007 to 2012. 
Dr. Sun is an editor of IEEE Transactions on Smart Grid and an associate editor of IET Generation, Transmission and Distribu-tion. 
His research interests include power system dynamics, stability and control and complex systems.
\end{IEEEbiography}

\begin{IEEEbiography} [{\includegraphics[width=1in,height=1.25in,clip,keepaspectratio]{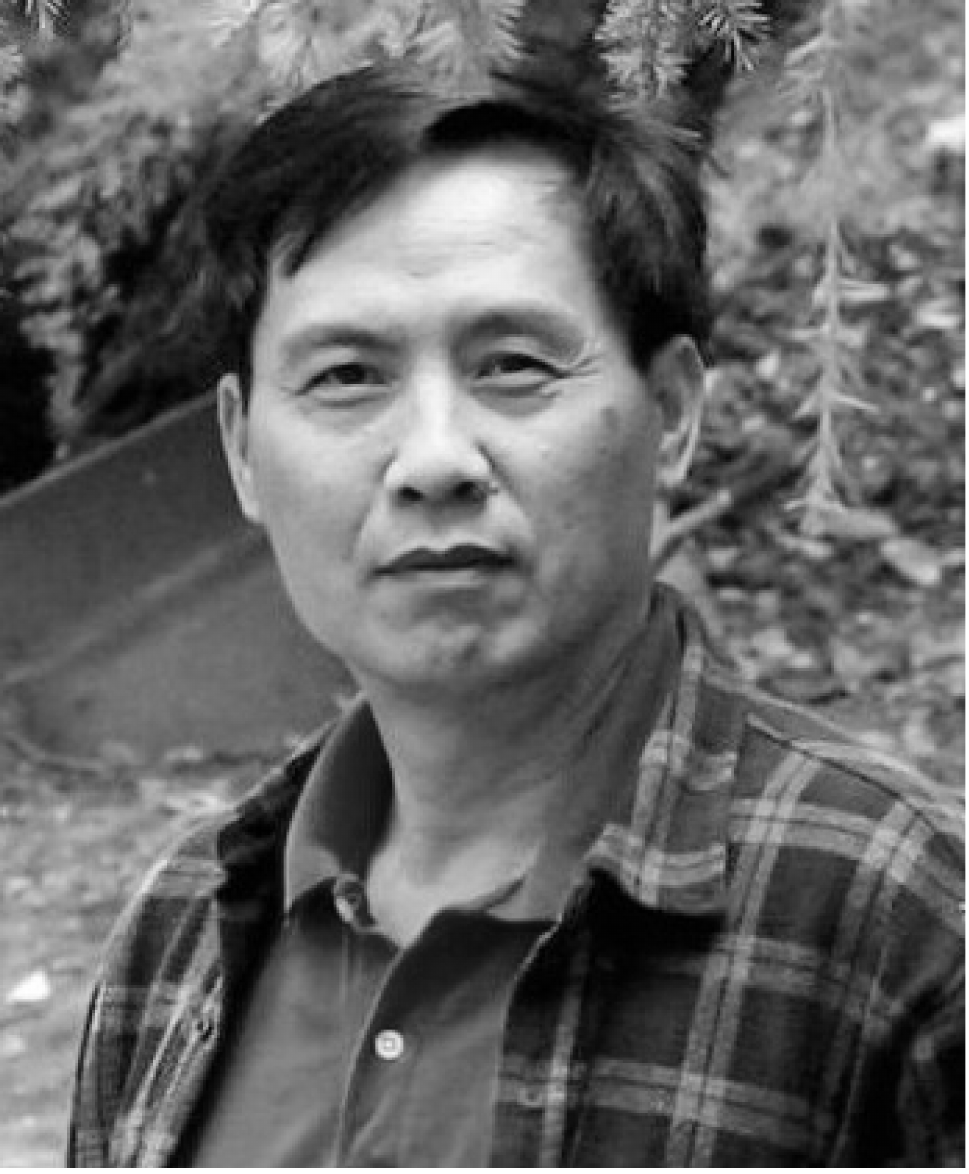}\vfill}]
{Wei Kang} (M'91--F'08)
received the B.S. and M.S. degrees from Nankai University, China, both in mathematics, in 1982 and 1985, respectively, and the Ph.D. degree in mathematics from the University of California, Davis, in 1991.

He is currently a Professor of applied mathematics at the U.S. Naval Postgraduate School, Monterey, CA.
He was a visiting Assistant Professor of systems science and mathematics at Washington University, St. Louis, MO (1991-1994).
He served as the Director of Business and International Collaborations at the American Institute of Mathematics (2008-2011).
His research interest includes computational optimal control, nonlinear filtering, cooperative control of autonomous vehicles, 
industry applications of control theory, nonlinear $H_\infty$ control, and bifurcations
and normal forms. His early research includes topics on Lie groups, Lie algebras, and differential geometry.

Dr. Kang is a fellow of IEEE. He was a plenary speaker in several international conferences of SIAM and IFAC. 
He served as an associate editor in several journals, including IEEE TAC and Automatica.
\end{IEEEbiography}

\end{document}